\documentclass[aos, preprint]{imsart}

\RequirePackage[OT1]{fontenc}
\RequirePackage{amsthm,amsmath}
\usepackage{array}
\usepackage{mathrsfs}
\usepackage{amssymb,bbm}
\usepackage{amsfonts}
\usepackage{amsmath,amssymb}
\usepackage{graphicx}
\usepackage{amscd}
\usepackage{epstopdf}
\usepackage{geometry}
\usepackage{color}
\usepackage{mathrsfs}
\usepackage{latexsym,bm}
\usepackage{array}
\usepackage{threeparttable}

\usepackage{threeparttable}
\usepackage{rotating}
\usepackage{booktabs}

\usepackage{rotating}
\usepackage{multirow}
\usepackage{amsmath}
\usepackage{subfigure}

\usepackage[plain,noend]{algorithm2e}

\DeclareFontFamily{U}{mathx}{\hyphenchar\font45}
\DeclareFontShape{U}{mathx}{m}{n}{
      <5> <6> <7> <8> <9> <10>
      <10.95> <12> <14.4> <17.28> <20.74> <24.88>
      mathx10
      }{}
\DeclareSymbolFont{mathx}{U}{mathx}{m}{n}
\DeclareMathAccent{\widecheck}{\mathalpha}{mathx}{"71}

\DeclareMathAccent{\widecheck}{\mathalpha}{mathx}{"71}

\startlocaldefs \numberwithin{equation}{section} \theoremstyle{it}
\newtheorem{thm}{Theorem}[section]

\newtheorem{ass}{Assumption}[section]

\newtheorem{rem}{Remark}
\endlocaldefs

\begin{document}

\begin{frontmatter}
\title{Time series models for realized covariance matrices based on the matrix-F distribution}
\begin{aug}
	\author{\fnms{Jiayuan} \snm{Zhou}\thanksref{m1}\ead[label=e1]{zhou.j@ufl.edu}},
	\author{\fnms{Feiyu} \snm{Jiang}\thanksref{m2}\ead[label=e2]{jfy16@mails.tsinghua.edu.cn}},
	\author{\fnms{Ke} \snm{Zhu}\thanksref{m3}
\ead[label=e3]{mazhuke@hku.hk}}
	\and
	\author{\fnms{Wai Keung} \snm{Li}\thanksref{m3}\thanksref{m4}
		\ead[label=e4]{waikeungli@eduhk.hk}
		\ead[label=u1,url]{http://www.foo.com}}

	
	\affiliation{University of Florida\thanksmark{m1}, Tsinghua University,\thanksmark{m2}\\The University of  Hong Kong\thanksmark{m3} and The Education University of Hong Kong\thanksmark{m4} }
	%
	%

	\address{University of Florida\\
		Department of Statistics\\
		Florida, U.S.A.\\
		\printead{e1}\\
	}
	
	\address{Center for Statistical Science\\
		\quad and Department of Industrial Engineering\\ Tsinghua University\\ Beijing 100084, China\\
		\printead{e2}\\
	}
	
	\address{The University of Hong Kong\\
		Department of Statistics \& Actuarial Science\\
		Hong Kong\\
		\printead{e3}\\}
	
	\address{The Education University of Hong Kong\\
	Department of Mathematics and Information Technology\\
	Hong Kong\\
	\printead{e4}\\}	
\end{aug}

\begin{abstract}
We propose a new \underline{C}onditional \underline{B}EKK matrix-\underline{F} (CBF) model for the time-varying realized covariance (RCOV) matrices. This CBF model is capable of capturing heavy-tailed RCOV, which is an important stylized fact but could not be handled adequately by
the Wishart-based models. To further mimic the long memory feature of the RCOV,
a special CBF model with the conditional heterogeneous autoregressive (HAR) structure is introduced.
Moreover, we give a systematical study on the probabilistic properties and statistical  inferences of the CBF model,  including
exploring its stationarity,
establishing the asymptotics of
its  maximum likelihood estimator, and giving some new inner-product-based tests for its model checking.
In order to handle a large dimensional RCOV matrix, we construct two reduced CBF models
--- the variance-target CBF model (for moderate but fixed dimensional RCOV matrix) and the factor CBF model
(for high dimensional RCOV matrix). For both reduced models, the asymptotic theory of the estimated parameters is derived.
The importance of our entire methodology is illustrated by simulation results and two real examples.
\end{abstract}



\end{frontmatter}

\newpage

\section{Introduction}

Modeling the multivariate volatility of many asset returns is crucial for asset pricing, portfolio selection, and risk management.
After the seminal work of Barndorff-Nielsen and Shephard (2002, 2004)
and Andersen et al. (2003), the realized covariance (RCOV) matrix, estimated
from the intra-day high frequency return data, has been recognized as a better estimator than the daily squared returns for daily volatility.
Consequently, increasing attention has been focused on the modeling and
forecasting of these RCOVs; see, e.g.,
McAleer and Medeiros (2008), Hansen et al. (2012), Noureldin et al. (2012), Bollerslev et al. (2016), and many others.

Existing models for the RCOV matrices can be roughly categorized into two types: transformation-based models and likelihood-based models.
Models in the first category capture the dynamics of
the RCOV matrices in an indirect way via transformation. Bauer and Vorkink (2011) used a factor model for the vectorization of the
log transformation of RCOV matrix;
Chiriac and Voev (2011) applied a vector autoregressive fractionally integrated moving average process to
model the Cholesky decomposition of RCOV matrix;
Callot et al. (2017) transformed the RCOV matrix into a large vector by the $vech$ operator, and then
fitted this transformed vector by a vector autoregressive model. In the first two models,
the dimension of RCOV matrix has to be moderate (e.g., less than 6) for a feasible manipulation.
In the third model, the dimension of RCOV matrix is allowed to be 30 in applications with the help of the LASSO method.

Models in the second category deals with RCOV matrices directly by assuming that the innovation, which
drives the RCOV time series, has a specific matrix distribution to generate random positive definite matrices automatically without imposing additional constraints. This important feature
results in positive-definite estimated RCOV matrices.
Unlike scalar or vector distributions, so far only few matrix distributions have been found to have explicit forms.
The primary choice for the innovation distribution is Wishart, leading to
the Wishart autoregressive (WAR) model  in Gouri\'{e}roux et al. (2009), the conditional autoregressive Wishart (CAW) model
in Golosnoy et al. (2012), the mixture Wishart model in Jin and Maheu (2013, 2016), and the generalized CAW model in Yu et al. (2017) to name a few.
The other choice for the innovation distribution is matrix-F, which was recently adopted by Opschoor et al. (2018). Generally speaking,
matrix-F distribution is the generalization of the usual F distribution, while Wishart distribution is the
generalization of the $\chi^2$ distribution (see, e.g., Konno (1991) and Opschoor et al. (2018) for more discussions). Therefore, matrix-F
distribution could be more appropriate than Wishart distribution in capturing the heavy-tailed innovation, which is an important
stylized fact in many applications (see, e.g., Bollerslev (1987), Fan et al.
(2014), Zhu and Li (2015), and Oh and Patton (2017)).
These likelihood models have at least three edges over the transformation-based models.
First, the likelihood-based models  will preserve the useful and important matrix structural information, which makes  them more interpretable compared with transformation-based models.
Second,  the number of estimated parameters in the transformation-based models has order $O(n^4)$, while the one in the likelihood-based models
has order $O(n^2)$, where $n$ is the dimension of the RCOV matrix. When $n$ is large, the likelihood-based models can bring more convenience
and a less daunting task in computation.  Third,
the likelihood-based models  make use of the likelihood function of the RCOV matrices, and hence
their statistical inference methods could be easily provided.

This paper contributes to the literature from three aspects.
First, we propose a new \underline{C}onditional \underline{B}EKK matrix-\underline{F} (CBF) model to
study the time-varying RCOV matrices. Our CBF model has matrix-F distributed innovations
with two degrees of freedom parameters $\nu_1$ and $\nu_2$. When $\nu_2\to\infty$, our CBF model reduces to the
CAW model (Golosnoy et al., 2012), which has Wishart distributed innovations.
Hence, the degrees of freedom $\nu_2$ is designed to capture the heavy-tailedness of the RCOV.
Since the RCOV is also well documented to have long memory phenomenon, we further introduce a
special CBF model which has a similar conditional heterogeneous autoregressive (HAR) structure as in Corsi (2009).
This special model is coined the CBF-HAR model. Although the CBF-HAR model is not formally a long memory model,
it gives rise to persistence in the RCOV time series.
Two real examples demonstrate that our CBF model (especially the CBF-HAR model) can have a significantly better forecasting performance than the corresponding CAW model, and hence a simple incorporation of $\nu_2$ to capture the heavy-tailed RCOV is necessary from a practical viewpoint.

Second, we provide a systematically statistical inference procedure for the CBF model. Specifically,
we explore its stationarity conditions,
establish the strong consistency and asymptotic normality of its maximum likelihood
estimator (MLE), and investigate some new inner-product-based tests for model diagnostic checking. Moreover, the performance of  our entire methodology is assessed by simulation studies.
Compared to the existing BEKK-type multivariate time series models, our proofs of the entire inference procedure are much involved, since
the CBF model is tailored for matrix time series. Particularly, our inner-product-based tests
seem to be  the first diagnostic checking tool for matrix time series models, and the related idea can be easily extended to other models.

Third, we construct two reduced CBF models --- the variance targeted (VT) CBF (VT-CBF) model and the factor CBF (F-CBF) model,
to handle moderately large and high dimensional RCOV matrix respectively. For both reduced models, the asymptotic theory of the estimated parameters is derived.
The dimension of the RCOV matrix is allowed to be a moderate but fixed number in the
VT-CBF model, while it is allowed to grow with the sample size $T$ and the intra-day sample size in the F-CBF model.
Therefore, this
makes the prediction of large dimensional RCOV matrices feasible in many cases.
The importance of both reduced models is illustrated by two real applications.

The remainder of this paper is organized as follows. Section 2 introduces the CBF model and
studies its  probabilistic properties. Section 3 investigates the asymptotics of the MLE. Section
4 presents inner-product-based tests to check the model adequacy. Two reduced CBF models and their
related asymptotic theories are provided in Section 5.
Some simulation studies are carried out in Section 6. Applications are given in Section 7. Section 8
concludes this paper.
Proofs of all theorems are relegated to the supplementary material.


Some notations are used throughout the paper. $I_{n}$ is the identity matrix of order $n$, and $\otimes$
represents the Kronecker product.    For an $n\times n$ matrix $A$, $tr(A)$ is its trace, $A'$ is its transpose,
$\left|A\right|$ is its determinant, $\rho(A)$ is its biggest eigenvalue, $\left\Vert A\right\Vert =\sqrt{tr(A'A)}$ is its Euclidean (or Frobenius) norm,
$\left\Vert A\right\Vert _{spec}=\sqrt{\rho(A'A)}$ is its spectral norm,  $vec(A)$
is a vector obtained by stacking all the columns of $A$,  $vech(A)$
is a vector obtained by stacking all columns of the lower triagular part of $A$,  and
$A^{\otimes2}=A\otimes A$.

\section{Model and Properties}

\subsection{Model specification}

Let $Y_{t}^*$ be the integrated volatility matrix of $n$ asset returns $X_{t}$ at time $t=1,...,T$.
After the seminal work of Barndorff-Nielsen and Shephard (2002, 2004) and Andersen et al. (2003),
the $n\times n$ positive definite realized covariance (RCOV) matrix $Y_{t}$ calculated from the high-frequency return data of $X_{t}$ has been widely applied to estimate $Y_{t}^*$ in the literature; see, e.g.,
Barndorff-Nielsen et al. (2011),
Lunde et al. (2016), A\"{i}t-Sahalia and Xiu (2017), Kim et al. (2018) and references therein.
Moreover, $Y_{t}$ is often viewed as a precise estimate for the conditional variances and covariances of these $n$ low-frequency asset returns $X_{t}$, and hence
how to predict $Y_{t}$ by some dynamic models is important in practice.
Motivated by this, a new dynamic model for $Y_t$ is proposed in the current paper.

Let $\mathcal{G}_{t}=\sigma(Y_{s}; s\leq t)$
be a filtration up to time $t$.
We assume that
\begin{flalign}\label{f_model}
Y_{t}=\Sigma_{t}^{1/2}\Delta_{t}\Sigma_{t}^{1/2},
\end{flalign}
where $\{\Delta_{t}\}_{t=1}^{T}$ is a sequence of independent and identically distributed (i.i.d.) $n\times n$ positive definite random innovation matrices with $E(\Delta_{t}|\mathcal{G}_{t-1})=I_{n}$, each $\Delta_{t}$ follows
the matrix-F distribution $F(\nu,\frac{\nu_{2}-n-1}{\nu_{1}}I_{n})$, and the density of $F(\nu,\Sigma)$ is
\begin{flalign}\label{2.2}
f(x;\nu,\Sigma)=
\Lambda(\nu)\times\frac{\left|\Sigma\right|^{-\nu_{1}/2}
\left|x\right|^{(\nu_{1}-n-1)/2}}{\left|I_{n}+\Sigma^{-1}x\right|^{(\nu_{1}+\nu_{2})/2}},\,\,\mbox{ for }x\in\mathcal{R}^{n\times n},
\end{flalign}
where $\nu=\left(\nu_{1},\nu_{2}\right)'$ with degrees of freedom
$\nu_{1}>n+1$ and $\nu_{2}>n+1$, $\Sigma$ is an $n\times n$ positive definite matrix, and
\begin{flalign*}
\Lambda(\nu)=\frac{\Gamma_{n}((\nu_{1}+\nu_{2})/2)}{\Gamma_{n}(\nu_{1}/2)\Gamma(\nu_{2}/2)}\,\, \mbox{ with }\,\,\Gamma_{n}(x)=\pi^{n(n-1)/4}\prod_{i=1}^{n}\Gamma(x+(1-i)/2);
\end{flalign*}
moreover, $\Sigma_{t}^{1/2}\in\mathcal{G}_{t-1}$ is the square root of the $n\times n$ positive definite matrix $\Sigma_{t}$,
which has a BEKK-type dynamic structure (see Engle and Kroner, 1995):
\begin{flalign}\label{bekk}
\Sigma_{t}=\Omega+\sum_{i=1}^{P}\sum_{k=1}^{K}A_{ki}Y_{t-i}A_{ki}'+\sum_{j=1}^{Q}\sum_{k=1}^{K}B_{kj}\Sigma_{t-j}B_{kj}',
\end{flalign}
where $\Omega$, $A_{ki}$, $B_{kj}$ are all $n\times n$
real matrices, the integers $P,Q,K$ are known as the orders of the model, and
$\Omega$ as well as the initial states
$\Sigma_{0},\Sigma_{-1},...,\Sigma_{-Q+1}$ are all positive
definite. Under model (\ref{f_model}),
\begin{flalign}\label{2.5}
Y_{t}|\mathcal{G}_{t-1}\sim F\left(\nu,\frac{\nu_{2}-n-1}{\nu_{1}}\Sigma_{t}\right)
\end{flalign}
with $E(Y_{t}|\mathcal{G}_{t-1})=\Sigma_{t}$, that is, the conditional distribution of $Y_{t}$ is matrix-F with a BEKK-type mean structure.
In this sense, we call model (\ref{f_model}) the \underline{C}onditional \underline{B}EKK matrix-\underline{F} (CBF) model.

The CBF model is related to the CAW model in Golosnoy et al. (2012), in which $\Delta_{t}$ follows the Wishart distribution.
 To see it clearly,
we follow Konno (1991) and Leung and Lo
(1996) to re-write $Y_{t}$ in model (\ref{f_model}) as
\begin{flalign}\label{decomposition}
Y_{t}=\left(\frac{\nu_{2}-n-1}{\nu_{1}}\right)\Sigma_{t}^{1/2}L_{t}^{1/2}R_{t}^{-1}L_{t}^{1/2}\Sigma_{t}^{1/2},
\end{flalign}
where $L_{t}\sim\mbox{Wishart}(\nu_{1},I_{n})$ and $R_{t}\sim\mbox{Wishart}(\nu_{2},I_{n})$
are independent. As
$\lim \limits_{\nu_{2}\to\infty}\nu_{2}^{-1}R_{t}$ $=I_{n}$ in probability, the identity (\ref{decomposition}) implies that when $\nu_{2}\rightarrow\infty$,
$Y_{t}|\mathcal{G}_{t-1}\sim \mbox{Wishart}(\nu_{1}, $ $\nu_{1}^{-1}\Sigma_{t})$,
which is  exactly the CAW model. Therefore, compared to the CAW model, the degrees of freedom $\nu_{2}$ in the CBF model accommodates the heavy-tailed RCOV, meaning that each $Y_{t,ij}$ from $Y_t$ satisfying (\ref{2.5}) could have a heavier tail than that from $Y_t$ satisfying $Y_{t}|\mathcal{G}_{t-1}\sim \mbox{Wishart}(\nu_{1}, \nu_{1}^{-1}\Sigma_{t})$ (see, e.g., Opschoor et al. (2018) for more discussions and examples). Clearly, the identity (\ref{decomposition}) also guarantees
$Y_{t}$ to be symmetric and positive definite, and it can
be used to generate $Y_{t}$ by using Wishart random variables.

\begin{figure}[htp]
    \includegraphics[width=35pc,height=20pc]{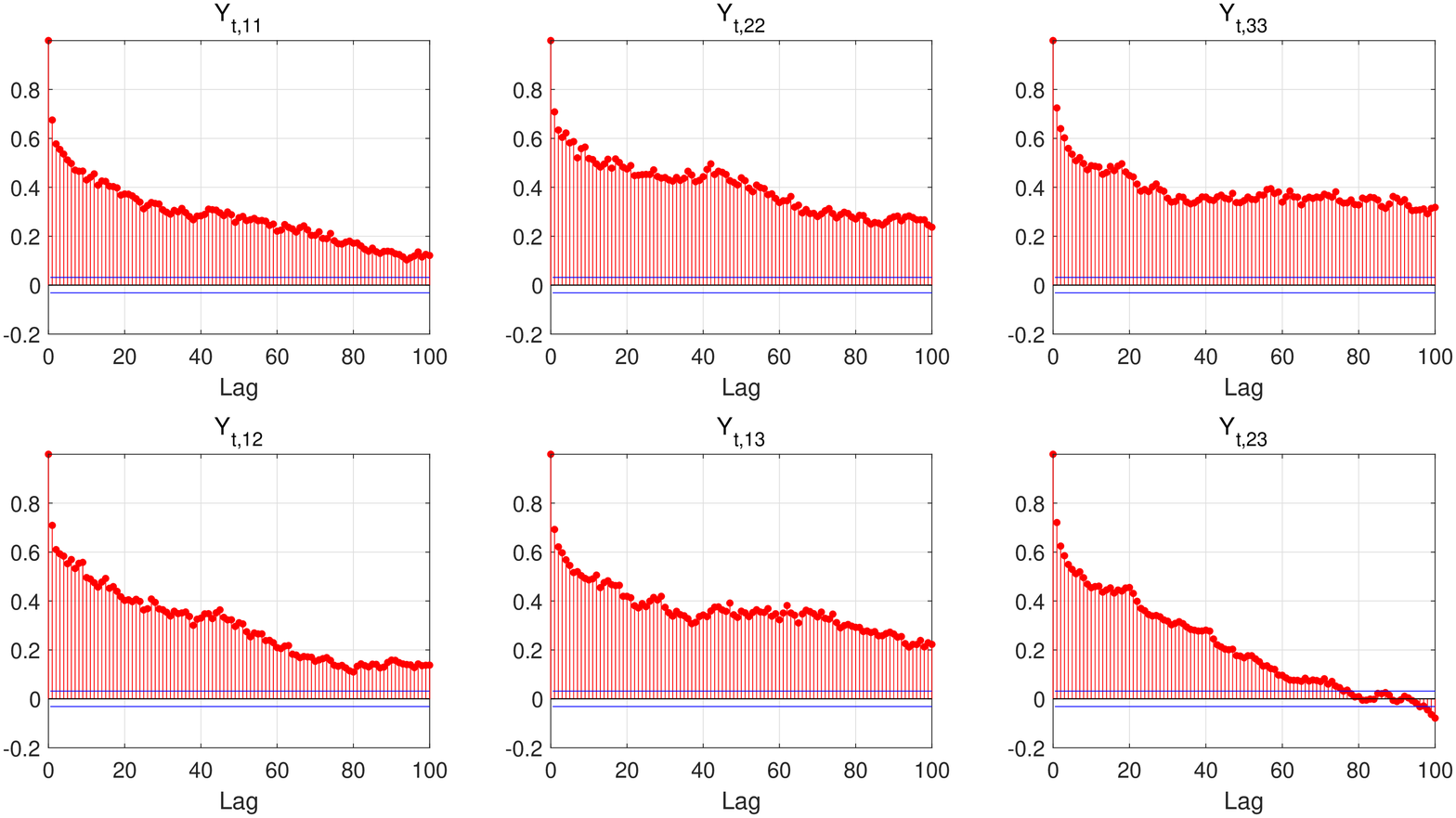}
    \caption{\label{fig0_1} Sample ACFs of one simulated data from a $3\times 3$ CBF-HAR model}
\end{figure}

Besides the heavy-tailedness, long memory is another well documented feature for the RCOV, and
it has been taken into account by many RCOV models, including
the heterogeneous autoregressive (HAR) model in Corsi (2009) as a benchmark. Although the HAR model does not formally belong to the class of
long memory models, it is able to reproduce the persistence of RCOV observed in many empirical data.
Inspired by the HAR model, we consider a special CBF model, which has the following specification for $\Sigma_t$:
\begin{flalign}\label{har_bekk}
\Sigma_{t}=\Omega+A_{(d)}Y_{t-1,d}A_{(d)}'+A_{(w)}Y_{t-1,w}A_{(w)}'+A_{(m)}Y_{t-1,m}A_{(m)}',
\end{flalign}
where $Y_{t-1,d}=Y_{t-1}$, $Y_{t-1,w}=(1/5)\sum_{i=1}^{5}Y_{t-i}$, and $Y_{t-1,m}=(1/22)\sum_{i=1}^{22}Y_{t-i}$ are
the daily, weekly, and monthly averages of RCOV matrices, respectively. In this case,
we label model (\ref{f_model}) as the CBF-HAR model, since
we put ``HAR dynamics'' on $\Sigma_t$.
Clearly, the CBF-HAR model is simply a constrained CBF model with $P=22$, $K=3$ and $Q=0$.
Figure \ref{fig0_1} plots the sample autocorrelation functions (ACFs) up to lag 100 of one simulated data from the  CBF-HAR model with
$\nu=(20,10)$ and
\begin{flalign*}
\Omega&=\left(\begin{array}{ccc}
0.5 & 0.2 &0.3 \\
0.2 & 0.5 &0.25\\
0.3& 0.25 &0.5
\end{array}\right), \quad\quad\,
A_{(d)}=\left(\begin{array}{ccc}
0.7 & 0 &0 \\
0 & 0.65&0\\
0&0&0.75
\end{array}\right),\\
A_{(w)}&=\left(\begin{array}{ccc}
0.6 & 0 &0 \\
0 & 0.6&0\\
0&0&0.55
\end{array}\right), \quad\quad
A_{(m)}=\left(\begin{array}{ccc}
0.4 & 0 &0 \\
0 & 0.45&0\\
0&0&0.4
\end{array}\right).
\end{flalign*}
From this figure, we can find that all entries of $Y_{t}$  exhibit
long memory phenomenon as expected.

Note that when $K=1$, sufficient identifiability conditions of model (\ref{bekk}) are that the main diagonal elements of $\Omega$ and
the first diagonal element of each $A_{1i}$, $B_{1j}$ are positive;
when $K>1$, some sufficient identifiability conditions of model (\ref{bekk}) can be found in Engle
and Kroner (1995).
For simplicity, we assume subsequently that model (\ref{bekk}) is identifiable.

Of course, the BEKK specification in model (\ref{bekk}) is not the only way to describe the dynamics of $\Sigma_{t}$. The
multivariate ARCH-type models such as the VEC model in Bollerslev
et al. (1988), the component model in Engle and Lee (1999),
the dynamic conditional correlation model in Engle (2002) and many others can also be adopted to model $\Sigma_{t}$.
Using these models together with the matrix-F distribution to fit and predict the RCOV matrices could be a promising direction for future study.

\subsection{Stationarity}

Stationarity is an important issue for most RCOV models, but so far it has been rarely studied. Denote $M=max(P,Q)$. For $i=1,2,\cdots,M$, let
$$A_{i}^{*}=\sum_{k=1}^{K}A_{ik}^{\otimes2}\,\,\,\mbox{ and }\,\,\,B_{i}^{*}=\sum_{k=1}^{K}B_{ik}^{\otimes2},$$
where $A_{ik}=0$ for $i>P$ and $B_{ik}=0$ for $i>Q$.
A sufficient condition for the stationarity of the CBF model is given below, and
it works for other general distributions of $\Delta_{t}$.

\begin{thm}\label{thm2.1}
Suppose that $\{\Delta_{t}\}$ in model (\ref{f_model}) is a sequence of i.i.d. $n\times n$ positive definite random matrices with $E\|\Delta_{t}\|<\infty$, and

\noindent (H1) the distribution of $\Delta_{1}$, denoted by $\varGamma$, is
absolute continuous with respect to the Lebesgue measure;

\noindent (H2) the point $I_{n}$ is in the interior of the support
of $\varGamma$;

\noindent  (H3) $\rho\left(\sum_{i=1}^{M}(A_{i}^{*}+B_{i}^{*})\right)<1$.

\noindent Then, $Y_{t}$ in model (\ref{f_model}) is strict stationary with $E\|Y_{t}\|<\infty$. Moreover,
$Y_{t}$ is positive Harris recurrent and geometrically ergodic.
\end{thm}

\begin{rem}
The results of Theorem \ref{thm2.1} are similar to those in Boussama
et al. (2011), where the stationarity of the BEKK model is studied.
Like Boussama
et al. (2011), the proof of Theorem \ref{thm2.1} is based on the semi-polynomial Markov chains technique, however, it
is much involved due to the matrix nature of model  (\ref{f_model}).
\end{rem}

As a special case,
the results in Theorem \ref{thm2.1} hold for the CAW model, in which $\Delta_{t}$ follows the Wishart distribution.
Under conditions (H1) and (H2), condition (H3) is  necessary and sufficient for the
strict stationarity of $Y_{t}$ with $E\|Y_{t}\|<\infty$. However, the
necessary and sufficient condition for the higher moments of $Y_t$ is still unclear at this stage.
Let $K_{n^2}$ be the $n^{2}\times n^{2}$ permutation matrix such that $K_{n^2}vec(A)=vec(A')$ for any $n\times n$ matrix $A$.
If $E\|Y_{t}\|^2<\infty$, it is not hard to see

\noindent (i) $\overline{y}:=E\left(vec(Y_{t})\right)=\Big[I_{n^{2}}-\sum_{i=1}^{M}\left(A_{i}^{*}+B_{i}^{*}\right)\Big]^{-1}vec(\Omega)$;

\noindent (ii)
$vec\left[E\left(vec(Y_{t})vec(Y_{t})'\right)\right]=\left(\Pi+I_{n^{4}}\right)
\left(I_{n^{4}}-\sum_{i=1}^{\infty}\Phi_{i}^{\otimes2}\Pi\right)^{-1}vec(\overline{y})\otimes vec(\overline{y})$,

\noindent where $\Pi=\left[s_{1}(\nu)-1\right]I_{n^{4}}+\left[s_{2}(\nu)I_{n^{2}}\otimes(I_{n^{2}}+K_{n^2})\right]\left[I_{n}\otimes K_{n^2}\otimes I_{n}\right]$ with
$$
s_{1}(\nu)=\frac{(\nu_{2}-n-1)[\nu_{1}(\nu_{2}-n-2)+2]}{\nu_{1}(\nu_{2}-n)(\nu_{2}-n-3)},\,\,\,
s_{2}(\nu)=\frac{(\nu_{2}-n-1)(\nu_{1}+\nu_{2}-n-1)}{\nu_{1}(\nu_{2}-n)(\nu_{2}-n-3)},
$$
and  $\Phi_{0}=I_{n^{2}}$, $\Phi_{i}=-B_{i}^{*}+\sum_{j=1}^{i}\big(A_{j}^{*}+B_{j}^{*}\big)\Phi_{i-j}$
for $i>0$. Result (ii) above clearly indicates that the parameters $\nu_1$ and $\nu_2$ have impact
on the second moment of $Y_t$ in a non-linear way. Although a closed form of third moment of $Y_t$ is absent,
similar impact from $\nu_1$ and $\nu_2$ is expected for the third moment of $Y_t$ and hence
 the asymptotic distribution of the proposed estimator (see Theorem \ref{t2} below).

\section{Maximum Likelihood Estimation}

Let $\theta=(\gamma',\nu')'\in\Theta$ be the unknown parameter of model (\ref{f_model}) with the true value $\theta_{0}=(\gamma_{0}',\nu_{0}')'$, where $\Theta=\Theta_{\gamma}\times\Theta_{\nu}$ is the parametric space with
$\Theta_{\gamma}\subset\mathbb{R}^{\tau_{1}}$ and $\Theta_{\nu}\subset\mathbb{R}^{2}$,
$\gamma=(w', u')'$, $w=vech(\Omega)$,
$u=(vec(A_{11})',...,vec(A_{KP})', vec(B_{11})',...,vec(B_{KQ})')$,  and $\tau_{1}=\frac{1}{2}n+[(P+Q)K+\frac{1}{2}]n^2$.
Below, we assume that $\Theta_{\gamma}$ and $\Theta_{\nu}$ are compact and $\theta_{0}$ is an interior point of $\Theta$.

Given the observations $\{Y_{t}\}_{t=1}^{T}$ and the initial values $\{Y_{t}\}_{t\leq 0}$,
the negative  log-likelihood function based on (\ref{2.5}) is
\begin{flalign}\label{3.1}
L(\theta)=\frac{1}{T}\sum_{t=1}^{T}l_{t}(\theta),
\end{flalign}
where
\begin{flalign*}
l_{t}(\theta)&=\frac{\nu_{1}}{2}\log\Big|\frac{\nu_2-n-1}{\nu_1}\Sigma_{t}(\gamma)\Big|-\frac{\nu_{1}-n-1}{2}\log |Y_{t}|\\
&\quad+\frac{\nu_{1}+\nu_{2}}{2}\log \Big|I_{n}+\frac{\nu_{1}}{\nu_{2}-n-1}\Sigma_{t}^{-1}(\gamma)Y_{t}\Big|+C(\nu)
\end{flalign*}
with $C(\nu)=-\log \Lambda(\nu)$  and $\Sigma_{t}(\gamma)$  calculated recursively by
\begin{flalign}\label{3.2}
\Sigma_{t}(\gamma)=
\Omega+\sum_{i=1}^{P}\sum_{k=1}^{K}A_{ki}Y_{t-i}A_{ki}'+\sum_{j=1}^{Q}\sum_{k=1}^{K}B_{kj}\Sigma_{t-j}(\gamma)B_{kj}'.
\end{flalign}
Clearly, $\Sigma_{t}(\gamma_{0})=\Sigma_{t}$.

As the initial values $\{Y_{t}\}_{t\leq 0}$ are not observable, we shall modify $L(\theta)$ as
\begin{flalign}\label{3.3}
\widehat{L}(\theta)=\frac{1}{T}\sum_{t=1}^{T}\widehat{l}_{t}(\theta),
\end{flalign}
where
$\widehat{l}_{t}(\theta)$ is defined in the same way as $l_{t}(\theta)$ with $\Sigma_{t}(\gamma)$ being
replaced by $\widehat{\Sigma}_{t}(\gamma)$, and
$\widehat{\Sigma}_{t}(\gamma)$ is calculated in the same way as $\Sigma_{t}(\gamma)$
based on a sequence of given constant matrices  $h:=\{Y_{0},\cdots,Y_{-M+1},\Sigma_{0},...,\Sigma_{-M+1}\}$.
The minimizer, $\widehat{\theta}=(\widehat{\gamma}',\widehat{\nu}')'$,  of $\widehat{L}(\theta)$ on $\Theta$ is called the
maximum likelihood estimator (MLE) of $\theta_{0}$. That is,
\begin{flalign}\label{full_mle}
\widehat{\theta}= (\widehat{\gamma}',\widehat{\nu}')'=
 \underset{\theta\in\Theta}{\arg\min}\,\,
\widehat{L}(\theta).
\end{flalign}

To study the asymptotic properties of $\widehat{\theta}$, we need two assumptions below.

\begin{ass}\label{asm1}
$Y_{t}$ is strictly stationary and ergodic.
\end{ass}

\begin{ass}\label{asm2}
For $\gamma\in\Theta_{\gamma}$, if $\gamma\not=\gamma_{0}$, $\Sigma_{t}(\gamma)\not=\Sigma_{t}(\gamma_{0})$ almost surely (a.s.) for all $t$.
\end{ass}

\noindent Assumption \ref{asm1} is standard, and Assumption \ref{asm2} which is
in line with Comte and Lieberman
(2003) and Hafner and Preminger (2009) is the identification condition.
The following two theorems give the consistency and asymptotic normality
of $\widehat{\theta}$, respectively.

\begin{thm}\label{t1}
Suppose that Assumptions \ref{asm1}-\ref{asm2} hold and $E\|Y_{t}\|<\infty$. Then,
$\widehat{\theta}\xrightarrow{a.s.}\theta_{0}\mbox{ as }T\to\infty$.
\end{thm}

\begin{thm}\label{t2}
Suppose that Assumptions \ref{asm1}-\ref{asm2} hold, $E\|Y_{t}\|^{3}<\infty$, and
\begin{flalign}\label{invert_O}
\mathcal{O}=E\left(\frac{\partial^2 l_{t}(\theta_{0})}{\partial\theta\partial\theta'}\right)\mbox{ is invertible}.
\end{flalign}
Then, $\sqrt{T}(\widehat{\theta}-\theta_{0})\xrightarrow{d} N(0, \mathcal{O}^{-1})\mbox{ as }T\to\infty$.
\end{thm}

Based on the observations $\{Y_{t}\}_{t=1}^{T}$ and  a sequence of given constant matrices $h$,
we can use the analytic expression of $\partial^2 l_{t}(\theta) /( \partial\theta\partial\theta')$ (see Appendix D in the supplementary material)  to
estimate $\mathcal{O}$ by its sample counterpart. As the univariate ARCH-type models,
the coefficients in the main diagonal line of $\Omega$ are positive to ensure the positive definite of $\Sigma_t$. Hence,
the classical $t$ or Wald test, which is constructed by the estimate of $\mathcal{O}$, can not be used to detect whether their values are zeros; see Li et al. (2018) for more discussions on this context.

\section{Model Diagnostic Checking}\label{checking}

Diagnostic tests are crucial for model checking in multivariate time series analysis; see, e.g.,
Li and McLeod (1981), Ling and Li (1997), Tse (2002) and many others. However,
no attempt has been made for the stationary matrix time series.
In this section, we propose some new inner-product-based tests to check the adequacy of model (\ref{f_model}).

Let $\mathfrak{Z}_{t}(\gamma) = vec\big(\Sigma_{t}^{-1/2}(\gamma)Y_{t}\Sigma_{t}^{-1/2}(\gamma)-I_{n}\big)$ be the vectorized residual for a given $\gamma$, and
$
\mathbf{b}_{t,j}(\gamma)=\mathfrak{Z}_{t}'(\gamma)\mathfrak{Z}_{t-j}(\gamma)$
be the inner product of two vectorized residuals at lag $j$.
Then, we stack $\mathbf{b}_{t,j}(\gamma)$ up to lag $l$ to construct $\mathcal{V}_{l}(\gamma)$, where
$$
\mathcal{V}_{l}(\gamma)=\frac{1}{T}\sum_{t=l+1}^{T}\big(
\mathbf{b}_{t,1}\left(\gamma\right),
\mathbf{b}_{t,2}\left(\gamma\right),
\cdots,
\mathbf{b}_{t,l}\left(\gamma\right)\big)',
$$
and $l\geq1$ is a given integer. Our testing idea is motivated by the fact that if model (\ref{f_model}) is adequate, $\mathfrak{Z}_{t}(\gamma_{0})$ is
a sequence of i.i.d. random vectors with mean zero, and hence
the value of $\mathcal{V}_{l}(\widehat{\gamma})$ is expected to be  close to zero.
To implement our test,  we need study the asymptotic property of $\mathcal{V}_{l}(\widehat{\gamma})$
in the following theorem.

\begin{thm}\label{cat_inner}
	Suppose that Assumptions \ref{asm1}-\ref{asm2} hold, $E\|Y_{t}\|^{4}<\infty$, and (\ref{invert_O}) holds.
	Then, if model (\ref{f_model}) is correctly specified,
	$
	\sqrt{T}\mathcal{V}_{l}(\widehat{\gamma})\xrightarrow{d}N(
	0, \mathbf{V})
	$ as $T\to\infty$,
	where
	$\mathbf{V} =(
	I_{l}, \mathfrak{R}_1)\mathfrak{R}_2(
	I_{l}, \mathfrak{R}_1)'
	$ with
\begin{flalign*}
	&\mathfrak{R}_1=E\left(\begin{array}{c}
	\mathfrak{Z}_{t-1}'\left(\gamma_{0}\right)\left(\partial\mathfrak{Z}_{t}\left(\gamma_{0}\right)/\partial\theta'\right)\\
	\mathfrak{Z}_{t-2}'\left(\gamma_{0}\right)\left(\partial\mathfrak{Z}_{t}\left(\gamma_{0}\right)/\partial\theta'\right)\\
	\vdots\\
	\mathfrak{Z}_{t-l}'\left(\gamma_{0}\right)\left(\partial\mathfrak{Z}_{t}\left(\gamma_{0}\right)/\partial\theta'\right)
	\end{array}\right) \times \mathcal{O}^{-1}
    \mbox{ and } \mathfrak{R}_2=\left(\begin{array}{cc}
tr\{E^2[\mathfrak{Z}_{t}'(\gamma_0)\mathfrak{Z}_{t}(\gamma_0)]\}I_l &0
\\0  & \mathcal{O}
\end{array}\right).
\end{flalign*}
\end{thm}


Based on Theorem \ref{cat_inner}, we construct the inner-product-based test statistic
\begin{flalign}\label{diag_test_inner}
\Pi(l)=T\big[\mathcal{V}_{l}'(\widehat{\gamma})\widehat{\mathbf{V}}^{-1}\mathcal{V}_{l}(\widehat{\gamma})\big]
\end{flalign}
to detect the adequacy of model (\ref{f_model}), where  $\widehat{\mathbf{V}}$ is the sample counterpart of
$\mathbf{V}$. If $\Pi(l)$
is larger than the upper-tailed critical value of $\chi^{2}(l)$, the fitted model (\ref{f_model})
is not adequate at a given significance level. Otherwise, it could be deemed as adequate.

Note that if we consider a test based on $\{\mathfrak{Z}_{t}(\widehat{\gamma})\}$ directly,
the resulting limiting distribution shall still be chi-squared, but its degrees of freedom increases fast with the dimension $n$.
To avoid this dilemma, we use the inner product of the residuals to construct our test $\Pi(l)$,
whose limiting distribution is independent of $n$.
This new idea is different from the portmanteau test in Ling and Li (1997) in which the test statistic  is constructed based on the auto-correlations of
the transformed scale residuals, while our test $\Pi(l)$ is based on the auto-covariances of the original vectorized
residuals. Clearly, our idea can be easily extended to the framework in Ling and Li (1997).	
Meanwhile, our inner-product-based test $\Pi(l)$
takes the auto-covariances of all entries of $\mathfrak{Z}_{t}(\widehat{\gamma})$ into account, while the idea of regression-based
test in Tse (2002) only considers one entry of $\mathfrak{Z}_{t}(\widehat{\gamma})$ at a time.
In view of this, we prefer to use the proposed inner-product idea for testing purpose.


\section{The Reduced CBF Models}

As the number of parameters in the CBF model is $O(n^2)$,
the estimation of the CBF model could be very computationally demanding when $n$ is large. This section introduces two
reduced CBF models, which are feasible in fitting RCOV matrices with a large $n$.

\subsection{The VT-CBF model}

This subsection proposes a reduced CBF model by
using the variance target (VT) technique in Engle and Mezrich (1996).
The idea of VT is to re-parameterize  the drift matrix $\Omega$ by using the theoretical mean of $Y_{t}$,
so that the estimation of $\Omega$ is excluded
in the implementation of the maximum likelihood estimation. Other related studies on the VT time series models can be
found in Francq et al. (2011) and Pedersen and Rahbek (2014).

To define our reduced model, we assume that $Y_{t}$ is strictly stationary with a finite mean $S=E(Y_{t})$.
By taking expectation on both sides of (\ref{bekk}), we have
\begin{flalign}\label{vt}
\Omega=S-\sum_{i=1}^{P}\sum_{k=1}^{K}A_{ki}SA_{ki}'-\sum_{j=1}^{Q}\sum_{k=1}^{K}B_{kj}SB_{kj}',
\end{flalign}
due to the fact that $S=E(Y_{t})=E(\Sigma_{t})$. With the help of (\ref{vt}), model (\ref{f_model}) becomes
\begin{flalign}\label{vt_model}
Y_{t}=\Sigma_{t}^{1/2}\Delta_{t}\Sigma_{t}^{1/2},
\end{flalign}
where all notations are inherited from model (\ref{f_model}), except that
\begin{flalign}\label{vt_bekk}
\Sigma_{t}&=S-\sum_{i=1}^{P}\sum_{k=1}^{K}A_{ki}SA_{ki}'-\sum_{j=1}^{Q}\sum_{k=1}^{K}B_{kj}SB_{kj}'\nonumber\\
&\quad+\sum_{i=1}^{P}\sum_{k=1}^{K}A_{ki}Y_{t-i}A_{ki}'+\sum_{j=1}^{Q}\sum_{k=1}^{K}B_{kj}\Sigma_{t-j}B_{kj}'.
\end{flalign}
We call model (\ref{vt_model}) the VT-CBF model.
Clearly, this reduced model shares the same probabilistic properties as the full CBF model. Although the VT-CBF model has the same amount of parameters as the full CBF model,
its two-step estimator given below is computationally easier
than the MLE for the full CBF model.

To present this two-step estimator, we let $\theta_{v}=(\delta',\nu')'\in \Theta_{v}$ be the unknown parameters of model (\ref{vt_model}) and its true value be $\theta_{v0}=(\delta_{0}',\nu_{0}')'$, where $\Theta_{v}=\Theta_{\delta}\times\Theta_{\nu}$ is
the parametric space with $\Theta_{\delta}=\Theta_{s}\times\Theta_{u}\subset\mathbb{R}^{\tau_{2}}$, $\tau_{2}=[(P+Q)K+1]n^2$ and $\Theta_{\nu}\subset\mathbb{R}^{2}$. Let
$\delta=(s',u')'$ with
$s= vec(S)$, $\Theta_s\in\mathbb{R}^{n^2}$ and $\Theta_u\in\mathbb{R}^{[(P+Q)Kn^2]}$.
As before, we assume that $\Theta_{\delta}$ and $\Theta_{\nu}$ are compact and $\theta_{v0}$ is an interior point of $\Theta_{v}$.

In the first step, we estimate
$s$ by $\widehat{s}_{v}$, where
$\widehat{s}_{v}=vec\left(\overline{Y_t}\right):=vec\big(\frac{1}{T}\sum_{t=1}^{T}Y_{t}\big)$.
In the second step, we estimate the remaining parameters $\zeta=(u',\nu')'$ by the constrained  MLE
based on the following modified log-likelihood function:
\begin{flalign}\label{3.5}
\widehat{L}_{v}(\theta_{v})=\frac{1}{T}\sum_{t=1}^{T}\widehat{l}_{vt}(\theta_{v}),
\end{flalign}
where
\begin{flalign*}
\widehat{l}_{vt}(\theta_{v})&=\frac{\nu_{1}}{2}\log\Big|\frac{\nu_2-n-1}{\nu_1}\widehat{\Sigma}_{vt}(\delta)\Big|-\frac{\nu_{1}-n-1}{2}\log |Y_{t}|\\
&\quad+\frac{\nu_{1}+\nu_{2}}{2}\log \Big|I_{n}+\frac{\nu_{1}}{\nu_{2}-n-1}\widehat{\Sigma}_{vt}^{-1}(\delta)Y_{t}\Big|+C(\nu),
\end{flalign*}
and $\widehat{\Sigma}_{vt}(\delta)$ is calculated recursively by
\begin{flalign}\label{recurs}
\widehat{\Sigma}_{vt}(\delta)&=S-\sum_{i=1}^{P}\sum_{k=1}^{K}A_{ki}SA_{ki}'-\sum_{j=1}^{Q}\sum_{k=1}^{K}B_{kj}SB_{kj}'\nonumber\\
&\quad+\sum_{i=1}^{P}\sum_{k=1}^{K}A_{ki}Y_{t-i}A_{ki}'+\sum_{j=1}^{Q}\sum_{k=1}^{K}B_{kj}\widehat{\Sigma}_{vt-j}(\delta)B_{kj}',
\end{flalign}
based on a sequence of given constant matrices  $h$.
Clearly, $\widehat{L}_{v}(\theta_{v})$ is analogous to  $\widehat{L}(\theta)$ in (\ref{3.3}),
and it is the modification of the following log-likelihood function:
\begin{flalign}\label{3.7}
L_{v}(\theta_{v})=\frac{1}{T}\sum_{t=1}^{T}l_{vt}(\theta_{v}),
\end{flalign}
where $l_{vt}(\theta_{v})$ is defined in the same way as $\widehat{l}_{vt}(\theta_{v})$ with $\widehat{\Sigma}_{vt}(\delta)$
being replaced by $\Sigma_{vt}(\delta)$,
and $\Sigma_{vt}(\delta)$ is calculated recursively by
\begin{flalign}\label{3.8}
\Sigma_{vt}(\delta)&=S-\sum_{i=1}^{P}\sum_{k=1}^{K}A_{ki}SA_{ki}'-\sum_{j=1}^{Q}\sum_{k=1}^{K}B_{kj}SB_{kj}'\nonumber\\
&\quad+\sum_{i=1}^{P}\sum_{k=1}^{K}A_{ki}Y_{t-i}A_{ki}'+\sum_{j=1}^{Q}\sum_{k=1}^{K}B_{kj}\Sigma_{vt-j}(\delta)B_{kj}',
\end{flalign}
based on the observations $\{Y_{t}\}_{t=1}^{T}$ and the initial values $\{Y_{t}\}_{t\leq 0}$.
The minimizer, $\widehat{\zeta}_{v}=(\widehat{u}_{v}',\widehat{\nu}_{v}')'$,  of $\widehat{L}_{v}(\widehat{s}_{v},\zeta)$ on $\Theta_{u}\times\Theta_{\nu}$ is the
constrained MLE of $(u_{0}', \nu_{0}')'$. That is,
\begin{flalign}\label{vt_mle}
 (\widehat{u}_{v}',\widehat{\nu}_{v}')'=
 \underset{\zeta\in\Theta_{u}\times\Theta_{\nu}}{\arg\min}\,\,
 \widehat{L}_{v}(\widehat{s}_{v},\zeta).
\end{flalign}

Now, we call $\widehat{\theta}_{v}=(\widehat{s}_{v}',\widehat{\zeta}_{v}')'$ the two-step estimator of
$\theta_{v}$ in model (\ref{vt_model}).
Let $\Psi(u)=\big(I_{n^{2}}-\sum_{i=1}^{M}A_{i}^{*}-\sum_{i=1}^{M}B_{i}^{*}\big)^{-1}\big(I_{n^{2}}-\sum_{i=1}^{M}B_{i}^{*}\big)$ and
$w_{t}(\theta_{v})=\Big(\begin{array}{c}
	\Psi(u) vec(Y_{t}-\Sigma_{vt}(\delta))\\
	\partial l_{vt}(\theta_{v})/\partial \zeta \end{array}\Big)$.
The following two theorems give the consistency and asymptotic normality
of $\widehat{\theta}_{v}$, respectively.

\begin{thm}\label{t4}
	Suppose that Assumptions \ref{asm1}-\ref{asm2} hold and $E\|Y_{t}\|<\infty$. Then,
	$\widehat{\theta}_{v}\xrightarrow{a.s.}\theta_{v0}\mbox{ as }T\to\infty$.
\end{thm}

\begin{thm}\label{t5}
	Suppose that Assumptions \ref{asm1}-\ref{asm2} hold, $E\|Y_{t}\|^{3}<\infty$, and
\begin{flalign}\label{invert_J}
J_{1}=E\Big[\frac{\partial^{2}l_{vt}(\theta_{v0})}{\partial \zeta\partial \zeta'}\Big]\mbox{ is invertible}.
\end{flalign}
Then,
	$\sqrt{T}(\widehat{\theta}_{v}-\theta_{v0})\xrightarrow{d} N(0, \mathcal{O}_{v})\mbox{ as }T\to\infty$,
where
\begin{flalign*}
	\mathcal{O}_{v}=\left(\begin{array}{cc}I_{n^2} & 0\\	-J_{1}^{-1}J_{2} & -J_{1}^{-1}\end{array}\right)	E(w_{t} w'_{t})
\left(\begin{array}{cc}	I_{n^2} & 0\\-J_{1}^{-1}J_{2} & -J_{1}^{-1}	\end{array}\right)'
	\end{flalign*}
with
	$J_{2}=E\Big[\frac{\partial^{2}l_{vt}(\theta_{v0})}{\partial \zeta\partial s'}\Big]$ and
$w_{t}=w_{t}(\theta_{v0})$.
\end{thm}

As before, we can use the sample counterpart of the analytic expressions of $\partial l_{vt}(\theta_{v})/\partial\theta_{v}$ and $\partial^2 l_{vt}(\theta_{v})/\partial\theta_{v}\partial\theta_{v}'$ to estimate
$\mathcal{O}_{v}$. Although the VT-CBF model can be estimated by the aforementioned two-step estimation procedure,
it still has to handle a large number of estimated parameters with order $O(n^{2})$ caused by the parameter matrices $A_{ki}$ and $B_{kj}$.
To make a more parsimonious  VT-CBF model, we can further impose some restrictions
on $A_{ki}$ and $B_{kj}$. McCurdy and Stengos (1992) and Engle and
Kroner (1995) have suggested to use
diagonal volatility models, which not only avoid over-parameterization, but also reflect
the fact that the variances and the covariances rely more on its own past
than the history of other variances or covariances. Motivated by this, we can assume that all $A_{ki}$ and $B_{kj}$ have a diagonal structure,
leading to a diagonal VT-CBF model. Clearly, the number of estimated parameters in the diagonal VT-CBF model has order $O(n)$, which is feasible to be handled for a moderate large but fixed $n$.

Next, similar to $\Pi(l)$ in (\ref{diag_test_inner}), we can construct the inner-product-based test statistics
to check the adequacy of model (\ref{f_model}) based on the two-step estimator $\widehat{\theta}_{v}$.
Let $\delta_{0}=(s_{0}',u_{0}')'$, $\widehat{\delta}_{v}=(\widehat{s}_{v}',\widehat{u}_{v}')'$,
$\mathfrak{Z}_{vt}(\delta) = vec(\Sigma_{vt}^{-1/2}(\delta)Y_{t}\Sigma_{vt}^{-1/2}(\delta)-I_{n})$ be the residual vector for a given $\delta$,
$
\mathbf{b}_{vt,j}(\delta)=\mathfrak{Z}_{vt}'(\delta)\mathfrak{Z}_{vt-j}(\delta)$
be the inner product of the residuals at lag $j$, and
$$
\mathcal{V}_{vl}(\delta)=\frac{1}{T}\sum_{t=l+1}^{T}\big(
\mathbf{b}_{vt,1}\left(\delta\right),
\mathbf{b}_{vt,2}\left(\delta\right),
\cdots,
\mathbf{b}_{vt,l}\left(\delta\right)\big)'.
$$
The asymptotic property of $\mathcal{V}_{vl}(\widehat{\delta}_{v})$ is given in the following theorem.

\begin{thm}\label{vt_cat_inner}
	Suppose that Assumptions \ref{asm1}-\ref{asm2} hold, $E\|Y_{t}\|^{4}<\infty$, and (\ref{invert_J}) holds. Then, if model
	(\ref{f_model}) is correctly specified, $\sqrt{T}\mathcal{V}_{vl}(\widehat{\delta}_{v})\xrightarrow{d}N(
	0, \mathbf{V}_{v})$ as $T\to\infty$, where
	$
	\mathbf{V}_{v}=(
	I_{l}, \mathfrak{R}_{1v})\mathfrak{R}_{2v}(
	I_{l}, \mathfrak{R}_{1v})'
	$ with
\begin{flalign*}
	&\mathfrak{R}_{1v}=E\left(\begin{array}{c}
	\mathfrak{Z}_{vt-1}'\left(\delta_{0}\right)\left(\partial\mathfrak{Z}_{vt}\left(\delta_{0}\right)/\partial\theta'\right)\\
	\mathfrak{Z}_{vt-2}'\left(\delta_{0}\right)\left(\partial\mathfrak{Z}_{vt}\left(\delta_{0}\right)/\partial\theta'\right)\\
	\vdots\\
	\mathfrak{Z}_{vt-l}'\left(\delta_{0}\right)\left(\partial\mathfrak{Z}_{vt}\left(\delta_{0}\right)/\partial\theta'\right)
	\end{array}\right)\times\left(\begin{array}{cc}
	I_{n^{2}} & 0\\
	-J_{1}^{-1}J_{2} & -J_{1}^{-1}
	\end{array}\right)
\end{flalign*}
and
\begin{flalign*}
	\mathfrak{R}_{2v}=\left(\begin{array}{cc}
		tr\{E^2[\mathfrak{Z}_{vt}(\delta_0)'\mathfrak{Z}_{vt}(\delta_0)]\}I_l &0
		\\0  &E(w_tw_t')
	\end{array}\right).
	\end{flalign*}
\end{thm}


By the preceding theorem, we can adopt the test statistic
\begin{flalign}\label{diag_test_inner_vt}
\Pi_{v}(l)=T[\mathcal{V}_{vl}'(\widehat{\delta}_{v})\widehat{\mathbf{V}}_{v}^{-1}\mathcal{V}_{vl}(\widehat{\delta}_{v})]
\end{flalign}
to detect the adequacy of model (\ref{f_model}), where  $\widehat{\mathbf{V}}_{v}$ is the sample counterpart of
$\mathbf{V}_{v}$. If $\Pi_{v}(l)$
is larger than the upper-tailed critical value of $\chi^{2}(l)$ at a given significance level, the fitted model (\ref{f_model})
is not adequate. Otherwise, it is adequate.

\subsection{The factor CBF model}

In modern data analysis, the dimension $n$ could be growing with the sample size $T$ in many cases, and this makes the CBF (or VT-CBF) models
computationally infeasible. Also, the dimension $n$ may be proportional to $m$ (the average intra-day sample size across all assets and all days), and then the methods to calculate $Y_{t}$ used for the fixed $n$ deliver an inconsistent estimator of $Y_{t}^{*}$; see, e.g., Wang and Zou (2010) and Tao et al. (2011) for surveys. To overcome this difficulty, we use the  thresholding
average realized volatility matrix (TARVM) estimator  in Tao et al. (2011) to calculate $Y_{t}$.
The TARVM is built based on the ARVM (Wang and Zou, 2010), which is constructed
 by taking the average of the constructed realized volatility matrices according to different predetermined  sampling frequencies.
The TARVM further thresholds the elements in each estimated RCOV matrix from the ARVM method,  so that certain sparsity structure is retained and the resulting estimator is consistent for
large $n$, which can be growing with (or even larger than) $T$.
For more recent works in this direction, we refer to  A\"{i}t-Sahalia and Xiu
(2017),  Kim et al. (2018),
and the references therein.

Since the dimension of $Y_{t}$ could be very large, it seems hard to study the dynamics of $Y_{t}$ without imposing some
specific structure. Here, we adopt the factor
model proposed by Tao et al. (2011) by assuming that
\begin{flalign}\label{factor_model}
Y_{t}^{*}=FY_{ft}^{*}F'+Y_{0}^{*},
\end{flalign}
where $Y_{ft}^{*}$ is an $r\times r$ positive definite factor covariance matrix with $r$ being a fixed integer (much smaller than $n$), $Y_{0}^{*}$ is an $n\times n$ positive definite constant matrix, and $F$ is an $n\times r$ factor loading matrix normalized by the constraint $F'F=I_{r}$. In model (\ref{factor_model}), the dynamic structure of
$Y_{t}^{*}$ is driven by that of a lower-dimensional latent process $Y_{ft}^{*}$, while $Y_{0}^{*}$
represents the static part of $Y_{t}^{*}$.

In (\ref{factor_model}), we shall highlight that only the column space of $F$ can be identified, and $F$ is not identified even if $F'F=I_r$ is imposed.
This is because $Y_{t}^{*}$ is unchanged when $F$ and $Y_{ft}^*$
are replaced by $F_{\dag}=FR$ and $Y_{ft,\dag}^{*}=R^{-1}{Y}_{ft}^*R^{-1'}$, respectively,
where $R$ is any $r\times r$ matrix satisfying $R'R=I_r$.


Define
\begin{flalign*}
\overline{Y}^{*}&=\frac{1}{T}\sum_{t=1}^{T}Y_{t}^{*},\,\,\,\,\,\overline{S}^{*}=\frac{1}{T}\sum_{t=1}^{T}\Big\{Y_{t}^{*}-\overline{Y}^{*}\Big\}^{2},
\end{flalign*}
and
\begin{flalign*}
\overline{Y}&=\frac{1}{T}\sum_{t=1}^{T}Y_{t},\,\,\,\,\,\overline{S}=\frac{1}{T}\sum_{t=1}^{T}\left\{Y_{t}-\overline{Y}\right\}^{2}.
\end{flalign*}
Then, we estimate $Y_{ft}^{*}$, $Y_{0}^{*}$ and $F$ by
\begin{flalign}\label{5.14}
\widehat{Y}_{ft}=\widehat{F}'Y_{t}\widehat{F},\,\,\,\,\widehat{Y}_{0}^{*}=\overline{Y}-\widehat{F}\widehat{F}'\overline{Y}\widehat{F}\widehat{F}'\,\,\,\mbox{ and }
\,\,\,\widehat{F}=(\widehat{f}_{1},\cdots,\widehat{f}_{r}),
\end{flalign}
respectively, where
$\widehat{f}_{1},\cdots,\widehat{f}_{r}$ are the eigenvectors of $\overline{S}$ corresponding to its
$r$ largest eigenvalues. As suggested by Lam and Yao (2012) and Ahn and Horenstein (2013), we may select $r$ such that the $r$ largest ratios of adjacent eigenvalues are significantly larger.


In order to study the asymptotics of the proposed estimators, we introduce the following technical assumptions.

\begin{ass}\label{asm3}
All row vectors of $F'$ and $Y_{0}^{*}$ satisfy the sparsity condition below. For an $n$-dimensional vector $(x_{1},\cdots,x_{n})$, we say
it is sparse if it satisfies
$$\sum_{i=1}^{n}|x_{i}|^{\delta_{*}}\leq U\pi(n),$$
where $\delta_{*}\in[0, 1)$, $U$ is a positive constant, and $\pi(n)$ is a deterministic function of $n$ that grows slowly in $n$
with typical examples $\pi(n)=1$ or $\log(n)$.
\end{ass}

\begin{ass}\label{asm4}
The factor model (\ref{factor_model}) has $r$ fixed factors, and matrices $Y_{0}^{*}$ and $Y_{ft}^{*}$ satisfy
$\|Y_{0}^{*}\|<\infty$ and
$
\max\limits_{1\leq t\leq T}\|Y_{ft,jj}^{*}\|=O_{p}(B(T))
$
for $j=1,2,\cdots,r$, where $Y_{ft,jj}^{*}$ is the $j$-th diagonal entry of $Y_{ft}^{*}$, and $1\leq B(T)=o(T)$.
\end{ass}

\begin{ass}\label{asm5}
$
\max\limits_{1\leq t\leq T}\|Y_{t}^{*}-Y_{t}\|=O_{p}(A(n,m,T))
$ for some rate function $A(n,m,T)$ such that $A(n,m,T)B^{5}(T)=o(1)$.
\end{ass}

Assumptions \ref{asm3}-\ref{asm5} are sufficient to prove the consistency of
$\widehat{Y}_{ft}$. For TARVM, we can take
$A(n,m,T)=\pi(n)[e_{m}(n^2T)^{1/\beta}]^{1-\delta_*}\log T$  and $B(T)=\log T$ with $e_m=m^{-1/6}$
so that $A(n,m,T)B^{5}(T)=o(1)$ for large $\beta$; see Tao et al. (2011).
Note that Assumptions \ref{asm3}-\ref{asm5} do not rule out the case that $n$ is larger than $T$, as long as $n^2T$ grows more slowly than $m^{\beta/6}$.
For other estimators, the rate $A(n,m,T)$ may be improved; see Tao et al. (2013) for more discussions.

\begin{thm}\label{t7}
Suppose that Assumptions \ref{asm3}-\ref{asm5} and the conditions in Theorem \ref{t2} hold.
Then, as $n, m, T$ go to infinity,

(i)\,\,$F'\widehat{F}-I_r=O_{p}(A(n,m,T)B(T))$,

(ii)\,\,$\widehat{Y}_{ft}-Y_{ft}=O_{p}(A^{1/2}(n,m,T)B^{3/2}(T))$,

\noindent where $Y_{ft}=Y_{ft}^{*}+F'Y_{0}^{*}F$, and $F=(f_{1},\cdots,f_{r})$ with
$f_{1},\cdots,f_{r}$ being the eigenvectors of $\bar{S}^{*}$ corresponding to its
$r$ largest eigenvalues.
\end{thm}

The above theorem indicates that $\widehat{Y}_{ft}$ is a consistent estimator of $Y_{ft}$ rather than $Y_{ft}^{*}$.
Next, we assume that $Y_{ft}$ satisfies the CBF model, that is,
\begin{flalign}\label{factor_cbf_star}
Y_{ft}|\mathcal{G}_{t-1}\sim F\left(\nu,\frac{\nu_{2}-n-1}{\nu_{1}}\Sigma_{ft}\right)
\end{flalign}
with $E(Y_{ft}|\mathcal{G}_{t-1})=\Sigma_{ft}$, where $\Sigma_{ft}$ is  defined in the same way as $\Sigma_{t}$ in (\ref{bekk}) with $Y_{t}$  replaced by $Y_{ft}$, and the remaining notations and set-ups inherent from model (\ref{f_model}).
We call models (\ref{factor_model}) and (\ref{factor_cbf_star}) the factor CBF (F-CBF) model.
Particularly, if $\Sigma_{ft}$ has the HAR dynamical structure as in (\ref{har_bekk}), the resulting model is called
the factor CBF-HAR (F-CBF-HAR) model. Based on this model, we have
$Y_{t}^{*}=F(Y_{ft}-F'Y_{0}^{*}F)F'+Y_{0}^{*}$. Since $Y_{t}\approx Y_{t}^{*}$, it implies that we
can study the large dimensional matrix $Y_{t}$ by using an $r\times r$ low-dimensional matrix $Y_{ft}$.

As $Y_{ft}$ is not observable,
we should estimate model (\ref{factor_cbf_star}) based on $\widehat{Y}_{ft}$, and hence we consider a feasible MLE of $\theta_{0}$ in model (\ref{factor_cbf_star}) given by
\begin{flalign*}
\widehat{\theta}_{1f}= (\widehat{\gamma}_{1f}',\widehat{\nu}_{1f}')'=
 \underset{\theta\in\Theta}{\arg\min}\,\,
\widehat{L}_{f}(\theta),
\end{flalign*}
where $\widehat{L}_{f}(\theta)$ is defined in the same way as $\widehat{L}(\theta)$ in (\ref{3.3}) with
$Y_{t}$ and $\widehat{\Sigma}_{t}(\gamma)$ replaced by $\widehat{Y}_{ft}$ and $\widehat{\Sigma}_{ft}(\gamma)$, respectively.
The following theorem shows that
$\widehat{\theta}_{1f}$ is consistent with the ideal MLE $\widehat{\theta}_{2f}$ based on $Y_{ft}$, where
\begin{flalign*}
\widehat{\theta}_{2f}= (\widehat{\gamma}_{2f}',\widehat{\nu}_{2f}')'=
 \underset{\theta\in\Theta}{\arg\min}\,\,
L_{f}(\theta),
\end{flalign*}
and $L_{f}(\theta)$ is defined in the same way as $L(\theta)$ in (\ref{3.1}) with $Y_{t}$ and $\Sigma_{t}(\gamma)$ replaced by $Y_{ft}$ and $\Sigma_{ft}(\gamma)$, respectively.

\begin{thm}\label{t8}
Suppose that the conditions in Theorem \ref{t7} hold.
Then, as $n, m, T$ go to infinity, $\widehat{\theta}_{1f}-\widehat{\theta}_{2f}=O_{p}(B(T)/T)+O_{p}(A^{1/2}(n,m,T)B^{5/2}(T))$.
\end{thm}

Since the dimension of $Y_{ft}$ is $r$ (much smaller than $n$), the calculation of $\widehat{\theta}_{1f}$  is computationally feasible.
In order to further reduce the number of parameters in model (\ref{factor_cbf_star}), we can also
assume that $Y_{ft}$ follows a VT-CBF model. This leads to the  F-VT-CBF model, which includes the F-VT-CBF-HAR model as a special case. For this F-VT-CBF model,
we consider its feasible two-step estimator $\widehat{\theta}_{1fv}=(\widehat{s}_{1fv}',\widehat{\zeta}_{1fv}')'$,  where
\begin{flalign*}
\widehat{s}_{1fv}=\frac{1}{T}\sum_{t=1}^{T}\widehat{Y}_{ft},\,\,\,\,\,
\widehat{\zeta}_{1fv}=(\widehat{u}_{1fv}',\widehat{\nu}_{1fv}')'=\underset{\zeta\in\Theta_{u}\times\Theta_{\nu}}{\arg\min}\,\,
 \widehat{L}_{fv}(\widehat{s}_{1fv},\zeta),
\end{flalign*}
and $\widehat{L}_{fv}(\theta_{v})$ is defined in the same way as $\widehat{L}_{v}(\theta_{v})$ in (\ref{3.5}) with
$Y_{t}$ and $\widehat{\Sigma}_{vt}(\delta)$ replaced by $\widehat{Y}_{ft}$ and $\widehat{\Sigma}_{fvt}(\delta)$, respectively.
Similar to Theorem \ref{t8},
$\widehat{\theta}_{1fv}$ is consistent with the ideal two-step estimator $\widehat{\theta}_{2fv}=(\widehat{s}_{2fv}',\widehat{\zeta}_{2fv}')'$ based on $Y_{ft}$, where
\begin{flalign*}
\widehat{s}_{2fv}=\frac{1}{T}\sum_{t=1}^{T}Y_{ft},\,\,\,\,\,
\widehat{\zeta}_{2fv}=(\widehat{u}_{2fv}',\widehat{\nu}_{2fv}')'=\underset{\zeta\in\Theta_{u}\times\Theta_{\nu}}{\arg\min}\,\,
 L_{fv}(\widehat{s}_{2fv},\zeta),
\end{flalign*}
and $L_{fv}(\theta_{v})$ is defined in the same way as $L(\theta_{v})$ in (\ref{3.7}) with $Y_{t}$ and $\Sigma_{t}(\delta)$ replaced by $Y_{ft}$ and $\Sigma_{fvt}(\delta)$, respectively.

\begin{thm}\label{t9}
Suppose that the conditions in Theorem \ref{t7} hold.
Then, as $n, m, T$ go to infinity,

(i)\,\,$\widehat{s}_{1fv}-\widehat{s}_{2fv}=O_{p}(A^{1/2}(n,m,T)B^{3/2}(T))$,

(ii)\,\,$\widehat{\zeta}_{1fv}-\widehat{\zeta}_{2fv}=O_{p}(B(T)/T)+O_{p}(A^{1/2}(n,m,T)B^{5/2}(T))$.
\end{thm}

Particularly, if $Y_{ft}$ follows a diagonal VT-CBF model,
the number of estimated parameters
in model (\ref{factor_cbf_star}) is $O(r)$, which is easy to calculate in practice.
In view of model (\ref{factor_model}) and the fact that $F'F=I_{r}$, we can predict $Y_{t}$ by either
$\widehat{F}\widehat{\Sigma}_{ft}(\widehat{\gamma}_{1f})\widehat{F}'+\widehat{Y}_{0}^{*}$ based on $\widehat{\theta}_{1f}$
or $\widehat{F}\widehat{\Sigma}_{fvt}(\widehat{\delta}_{1fv})\widehat{F}'+\widehat{Y}_{0}^{*}$ based on
$\widehat{\theta}_{1fv}$, where $\widehat{\delta}_{1fv}=(\widehat{s}_{1fv}',\widehat{u}_{1fv}')'$.

\section{Simulation}

In this section, we first assess the performance of the MLE $\widehat{\theta}$ and the two-step estimator $\widehat{\theta}_{v}$ in the finite sample.
We generate 1000 replications of sample size $T=1000$
and $2000$ from the following model:
\begin{flalign}\label{sim_model}
Y_{t}=\Sigma_{t}^{1/2}\Delta_{t}\Sigma_{t}^{1/2} \text{ with }\  \Sigma_{t}=\Omega_{0}+A_{10}Y_{t-1}A_{10}'+B_{10}\Sigma_{t-1}B_{10}',
\end{flalign}
where
\begin{flalign*} 
\Omega_0=\left(\begin{array}{ccc}
0.5 & 0.2 & 0.3\\
0.2 & 0.5 & 0.25\\
0.3 & 0.25 & 0.5
\end{array}\right),\,\, A_{10}=\left(\begin{array}{ccc}
0.4 & 0 & 0\\
0 & 0.55 & 0\\
0 & 0 & 0.5
\end{array}\right),\,\, B_{10}=\left(\begin{array}{ccc}
0.4 & 0 & 0\\
0 & 0.3 & 0\\
0 & 0 & 0.5
\end{array}\right),	
\end{flalign*}
$\{\Delta_{t}\}$ is a sequence of independent  $F\big(\nu_0,\frac{\nu_{20}-n-1}{\nu_{10}}I_{n}\big)$ distributed random matrices with $n=3$, and
$\nu_0=(10, 8), (15, 10)$ or $(20, 10)$. For each repetition, we calculate $\widehat{\theta}$, $\widehat{\theta}_{v}$, and their related asymptotic standard deviations. For $\widehat{\theta}_{v}$, we report the results related to $\Omega$ instead of $S$, and hence
the asymptotic standard deviation of the estimated parameters in $\Omega$ is absent in this case.

Table \ref{table1} reports the sample bias, the sample standard deviation (SD) and
the average asymptotic standard deviation (AD) of $\widehat{\theta}$ and $\widehat{\theta}_{v}$.
From this table, we can see that the biases of both estimators are small comparing to the magnitude of the parameters, and they become smaller as the sample size $T$ increases. This assures the accuracy of both estimators.
Furthermore, we find that the SDs are generally close to the ADs for both estimators, and all of the SDs and ADs become smaller as
$T$ increases from 1000 to 2000.
In terms of ADs or SDs, $\widehat{\theta}$ is generally more efficient than $\widehat{\theta}_{v}$, although this efficiency advantage is weak for many parameters. However, the estimation time for $\widehat{\theta}_{v}$ is almost 70\% of the time for $\widehat{\theta}$, and this computation advantage can be more significant when $n$ increases.

\begin{sidewaystable}
	\centering
	\caption{\label{table1} The results of the MLE $\widehat{\theta}$ and two-step estimator $\widehat{\theta}_{v}$ for model (\ref{sim_model})}
\addtolength{\tabcolsep}{-1.3pt}
\renewcommand{\arraystretch}{1}
	\scalebox{1}{
		\begin{tabular}{cc cc cc ccc ccc cccccc>{\centering}p{10cm}}
			\hline
		& T    &   &   & $\nu_1$ & $\nu_2$ & $A_{1,11}$ & $A_{1,22}$ & $A_{1,33}$ & $B_{1,11}$ & $B_{1,22}$ & $B_{1,33}$ & $\Omega_{11}$ & $\Omega_{21}$ & $\Omega_{31}$ & $\Omega_{22}$ & $\Omega_{32}$ & $\Omega_{33}$ \\ \hline
		Case 1 & 1000 &$\widehat{\theta}$& Bias & 0.0320        & 0.0160        & -0.0014        & -0.0029        & -0.0009        & -0.0151        & -0.0112        & -0.0102        & -0.0005             & 0.0028              & 0.0057              & -0.0009             & 0.0037              & 0.0053              \\
		&      && ESD  & 0.3914        & 0.2452        & 0.0255         & 0.0249         & 0.0240         & 0.1170         & 0.0964         & 0.0728         & 0.0600              & 0.0188              & 0.0337              & 0.0419              & 0.0248              & 0.0601              \\
		&      && ASD  & 0.4111        & 0.2563       & 0.0258         & 0.0259         & 0.0241         & 0.1103         & 0.0892         & 0.0652         & 0.0586              & 0.0179              & 0.0323              & 0.0402              & 0.0232              & 0.0562              \\
&  &$\widehat{\theta}_{v}$& Bias & -0.0080       & 0.0382        & -0.0005        & -0.0030        & 0.0000         & -0.0130        & -0.0088        & -0.0096        & -0.0020             & 0.0020              & 0.0047              & -0.0030             & 0.0033              & 0.0040              \\
		&      && ESD  & 0.3884        & 0.2607        & 0.0263         & 0.0272         & 0.0255         & 0.1165         & 0.0956         & 0.0728         & 0.0614              & 0.0229              & 0.0366              & 0.0433              & 0.0291              & 0.0615              \\
		&      && ASD  & 0.4024        & 0.2619        & 0.0266         & 0.0282         & 0.0258         & 0.1207         & 0.1046         & 0.0742         &                     &                     &                     &                     &                     &                     \\
		& 2000 &$\widehat{\theta}$&Bias & 0.0188        & 0.0072        & -0.0003        & -0.0018        & -0.0002        & -0.0111        & -0.0036        & -0.0034        & 0.0019              & 0.0014              & 0.0030              & -0.0008             & 0.0012              & 0.0014              \\
		&      &&ESD  & 0.2767        & 0.1733        & 0.0174         & 0.0179         & 0.0168         & 0.0797         & 0.0633         & 0.0459         & 0.0431              & 0.0130              & 0.0231              & 0.0293              & 0.0169              & 0.0405              \\
		&      &&ASD  & 0.2880        & 0.1797        & 0.0181         & 0.0182         & 0.0169         & 0.0767         & 0.0615         & 0.0447         & 0.0417              & 0.0124              & 0.0226              & 0.0287              & 0.0162              & 0.0395              \\
		&  &$\widehat{\theta}_{v}$& Bias & -0.0020       & 0.0196        & 0.0002         & -0.0020        & 0.0003         & -0.0103        & -0.0024        & -0.0030        & 0.0010              & 0.0007              & 0.0023              & -0.0020             & 0.0008              & 0.0005              \\
		&      && ESD  & 0.2767        & 0.1876        & 0.0181         & 0.0194         & 0.0179         & 0.0800         & 0.0633         & 0.0458         & 0.0438              & 0.0161              & 0.0250              & 0.0304              & 0.0200              & 0.0416              \\
		&      && ASD  & 0.2856        & 0.1871        & 0.0187         & 0.0198         & 0.0180         & 0.0823         & 0.0640         & 0.0458         &                     &                     &                     &                     &                     &                     \\
		&      &      &               &               &                &                &                &                &                &                &                     &                     &                     &                     &                     &                     \\
		Case 2 & 1000 &$\widehat{\theta}$& Bias & 0.0900        & 0.0132        & -0.0017        & -0.0021        & -0.0010        & -0.0159        & -0.0075        & -0.0092        & -0.0023             & 0.0026              & 0.0049              & -0.0022             & 0.0031              & 0.0053              \\
		&      && ESD  & 0.8099        & 0.3597        & 0.0240         & 0.0245         & 0.0227         & 0.1242         & 0.0974         & 0.0690         & 0.0626              & 0.0177              & 0.0340              & 0.0418              & 0.0237              & 0.0591              \\
		&      && ASD  & 0.8413        & 0.3598        & 0.0250         & 0.0243         & 0.0227         & 0.1158         & 0.0921         & 0.0672         & 0.0602              & 0.0175              & 0.0328              & 0.0410              & 0.0230              & 0.0578              \\
		 & &$\widehat{\theta}_{v}$& Bias & 0.0255        & 0.0353        & -0.0011        & -0.0019        & -0.0009        & -0.0154        & -0.0067        & -0.0089        & -0.0027             & 0.0025              & 0.0041              & -0.0031             & 0.0023              & 0.0041              \\
		&      && ESD  & 0.7985        & 0.3659        & 0.0244         & 0.0259         & 0.0232         & 0.1239         & 0.0974        & 0.0689         & 0.0632              & 0.0200              & 0.0352              & 0.0430              & 0.0257              & 0.0593              \\
		&      && ASD  & 0.8162        & 0.3585       & 0.0252         & 0.0253         & 0.0234         & 0.1290         & 0.1034         & 0.0724         &                     &                     &                     &                     &                     &                     \\
& 2000 &$\widehat{\theta}$& Bias & 0.0702        & 0.0000        & -0.0005        & -0.0011        & 0.0000         & -0.0074        & -0.0021        & -0.0041        & -0.0008             & 0.0009              & 0.0020              & -0.0017             & 0.0010              & 0.0020              \\
		&      && ESD  & 0.5912        & 0.2515        & 0.0173         & 0.0174         & 0.0158         & 0.0801         & 0.0665          & 0.0467         & 0.0434              & 0.0123              & 0.0228              & 0.0299              & 0.0163              & 0.0412              \\
		&      && ASD  & 0.5871        & 0.2517        & 0.0175         & 0.0171         & 0.0159         & 0.0805         & 0.0640         & 0.0461         & 0.0437              & 0.0121              & 0.0230              & 0.0293              & 0.0161              & 0.0406              \\
		&  &$\widehat{\theta}_{v}$& Bias & 0.0384        & 0.0112        & -0.0002        & -0.0009        & -0.0001        & -0.0071       & -0.0018        & -0.0040        & -0.0010             & -0.0007             & 0.0016              & -0.0021             & 0.0005              & 0.0012              \\
		&      && ESD  & 0.5874        & 0.2532        & 0.0175         & 0.0182         & 0.0163         & 0.0800         & 0.0664         & 0.0467         & 0.0437              & 0.0138              & 0.0235              & 0.0305              & 0.0176              & 0.0413              \\
		&      && ASD  & 0.5792       & 0.2537        & 0.0177         & 0.0178         & 0.0163         & 0.0827         & 0.0668         & 0.0472         &                     &                     &                     &                     &                     &                     \\
		&      &&      &               &               &                &                &                &                &                &                &                     &                     &                     &                     &                     &                     \\
		Case 3 & 1000 &$\widehat{\theta}$& Bias & 0.1521        & 0.0294        & -0.0013        & -0.0021        & -0.0007        & -0.0165        & -0.0092        & -0.0082        & -0.0023             & 0.0026              & 0.0046              & -0.0026             & 0.0032              & 0.0031              \\
		&      && ESD  & 1.4019        & 0.3340        & 0.0237         & 0.0237         & 0.0213         & 0.1253         & 0.0979         & 0.0712         & 0.0615              & 0.0173              & 0.0353              & 0.0418              & 0.0231              & 0.0599              \\
		&      && ASD  & 1.4496        & 0.3442        & 0.0242         & 0.0235         & 0.0220         & 0.1127         & 0.0904         & 0.0654         & 0.0586              & 0.0169              & 0.0319              & 0.0399              & 0.0224              & 0.0463              \\
&  &$\widehat{\theta}_{v}$& Bias & 0.0433        & 0.0475        & -0.0008        & -0.0013        & -0.0006        & -0.0156        & -0.0083        & -0.0082        & -0.0031             & 0.0025              & 0.0040              & -0.0031             & 0.0026              & 0.0021              \\
		&      && ESD  & 1.3774        & 0.3416        & 0.0239         & 0.0246         & 0.0220         & 0.1250         & 0.0975         & 0.0725         & 0.0614              & 0.0190              & 0.0358              & 0.0423              & 0.0250              & 0.0605              \\
		&      && ASD  & 1.4084        & 0.3443        & 0.0248         & 0.0246         & 0.0227         & 0.1409         & 0.0996         & 0.0755         &                     &                     &                     &                     &                     &                     \\
		& 2000 &$\widehat{\theta}$& Bias & 0.0737        & 0.0190        & -0.0006        & -0.0009        & 0.0001         & -0.0061        & -0.0047        & -0.0052        & -0.0016             & 0.0010              & 0.0018              & -0.0010             & 0.0016              & 0.0022              \\
		&      && ESD  & 1.0087        & 0.2469        & 0.0169         & 0.0163         & 0.0149         & 0.0794         & 0.0671         & 0.0480         & 0.0418              & 0.0118              & 0.0227              & 0.0295              & 0.0161              & 0.0418              \\
		&      && ASD  & 1.0057        & 0.2411        & 0.0170         & 0.0165         & 0.0155         & 0.0787         & 0.0630         & 0.0453         & 0.0429              & 0.0117              & 0.0225              & 0.0286              & 0.0157              & 0.0397              \\
&  &$\widehat{\theta}_{v}$& Bias & 0.0192        & 0.0286        & -0.0004        & -0.0004        & 0.0000         & -0.0058        & -0.0045        & -0.0051        & -0.0020             & 0.0008              & 0.0013              & -0.0012             & 0.0011              & 0.0015              \\
		&      && ESD  & 1.0022        & 0.2511        & 0.0170         & 0.0171         & 0.0153         & 0.0795         & 0.0673         & 0.0480         & 0.0419              & 0.0131              & 0.0231              & 0.0300              & 0.0176              & 0.0418              \\
		&      && ASD  & 0.9923        & 0.2434        & 0.0172         & 0.0173         & 0.0158         & 0.0817         & 0.0652         & 0.0462         &                     &                     &                     &                     &                     &\\
             \hline
	\end{tabular}}
  \begin{tablenotes}
      \item[\dag] {\scriptsize Cases 1-3 correspond to $\nu_0=(10, 8), (15, 10)$ and $(20, 10)$, respectively.}
  \end{tablenotes}
\end{sidewaystable}


Next, we examine the performance of the inner-product-based tests $\Pi(l)$ and $\Pi_{v}(l)$ in the finite sample.
We generate 1000 replications of sample size $T=1000$
and $2000$ from the following model:
\begin{flalign}\label{sim_model_2}
Y_{t}=\Sigma_{t}^{1/2}\Delta_{t}\Sigma_{t}^{1/2} \text{ with }  \Sigma_{t}=\Omega_0+A_{10}Y_{t-1}A_{10}'+A_{20}Y_{t-2}A_{20}'+B_{10}\Sigma_{t-1}B_{10}',
\end{flalign}
where the values of $\Omega_0$,
$A_{10}$ and $B_{10}$ are chosen as in (\ref{sim_model}),
$A_{20}=\mbox{diag}\{\lambda,\lambda,\lambda\}$ is a diagonal matrix with $\lambda=0, 0.05, 0.1, 0.15, 0.2$, and
$\{\Delta_{t}\}$ is a sequence of independent  $F\big(\nu_0,\frac{\nu_{20}-n-1}{\nu_{10}}I_{n}\big)$ distributed random matrices with $n=3$ and
$\nu_0=(10, 8)$. We fit each replication by the CBF model with $(K, P, Q)=(1, 1, 1)$, and use $\Pi(l)$ and $\Pi_{v}(l)$
to check whether the fitted model is adequate. Here, we set the significance level $\alpha=0.05$ and $l=2,3,4,5,6$. The empirical
sizes and powers of both tests are reported in Table \ref{innertest_size}, and
their sizes correspond to the results for the case of $\lambda=0$. From Table \ref{innertest_size}, we can find that both $\Pi(l)$ and $\Pi_{v}(l)$ always have accurate sizes, although they are slightly oversized for small $T$. For the power of both tests,
it is generally as expected. First, all the power becomes larger as $T$ increases. Second, both tests become more powerful as $\lambda$ becomes larger. Third, the power of $\Pi(l)$ and $\Pi_{v}(l)$ is comparable,
but the former need more computational time. Note that when $\nu_0=(15, 10)$ and $(20, 10)$, the testing results are similar to these
for $\nu_0=(10, 8)$, and hence they are not reported for saving space.

\begin{table}[!htb]
	\caption{\label{innertest_size}The results of $\Pi(l)$ and $\Pi_{v}(l)$ for model (\ref{sim_model_2})}
	\centering
\addtolength{\tabcolsep}{-1pt}
\renewcommand{\arraystretch}{1.2}
	\scalebox{1}{
\begin{tabular}{cc ccccc ccccc cccc}
	\hline
	& & \multicolumn{2}{c}{$l=2$} && \multicolumn{2}{c}{$l=3$} && \multicolumn{2}{c}{$l=4$} && \multicolumn{2}{c}{$l=5$} && \multicolumn{2}{c}{$l=6$}\\
\cline{3-4} \cline{6-7} \cline{9-10} \cline{12-13} \cline{15-16}
	$\lambda$ & $T$ & $\Pi(l)$ & $\Pi_{v}(l)$ && $\Pi(l)$ & $\Pi_{v}(l)$ && $\Pi(l)$ & $\Pi_{v}(l)$ && $\Pi(l)$ & $\Pi_{v}(l)$ && $\Pi(l)$ & $\Pi_{v}(l)$\\
\hline
	0 & 1000 & 0.043 & 0.037 && 0.048 & 0.045 && 0.052 & 0.054 && 0.047 & 0.048 && 0.049 & 0.054\\
	& 2000 & 0.048 & 0.056 && 0.058 & 0.059 && 0.053 & 0.054 && 0.052 & 0.059 && 0.051 & 0.052\\\\

	0.05 & 1000 & 0.048 & 0.045 && 0.051 & 0.048 && 0.058 & 0.053 && 0.060 & 0.052  && 0.061 & 0.062\\
	& 2000 & 0.060 & 0.063 && 0.063 & 0.073 && 0.064 & 0.075 && 0.063& 0.076 && 0.058 & 0.074\\\\

	0.1 & 1000  & 0.238 & 0.238 && 0.210 & 0.211 && 0.196 & 0.199 && 0.196 & 0.199 && 0.179 & 0.183 \\
	& 2000 & 0.414 & 0.408 && 0.371 &  0.364 && 0.350  &   0.354 && 0.309 & 0.328 && 0.316 & 0.320\\\\

	0.15 & 1000 & 0.885 & 0.854 && 0.847 & 0.818 && 0.818 & 0.793 && 0.784 &   0.762 && 0.768 & 0.746\\
	& 2000 &  0.974 & 0.956 && 0.966 & 0.951  && 0.956 & 0.933  && 0.946 & 0.925 && 0.941 & 0.919\\\\

	0.2 & 1000 & 0.976  & 0.924 && 0.972 & 0.916 && 0.964 & 0.893 && 0.961 & 0.889  && 0.956 & 0.887\\
	& 2000 & 0.992  & 0.951 && 0.989 & 0.945 && 0.987 & 0.923 && 0.987 & 0.914 && 0.985 & 0.910\\
    \hline
\end{tabular}}
\end{table}

Overall, both estimators $\widehat{\theta}$ and $\widehat{\theta}_{v}$  and both tests $\Pi(l)$ and $\Pi_{v}(l)$ have a good performance especially when the sample size $T$ gets larger.
When the dimension of $Y_{t}$ is small, our simulation results show that $\widehat{\theta}_{v}$ is only slightly less efficient than $\widehat{\theta}$, and $\Pi_{v}(l)$ is generally as powerful as $\Pi(l)$.
When the dimension of $Y_{t}$ is large, $\widehat{\theta}_{v}$ and $\Pi_{v}(l)$ can enjoy a faster computation speed than $\widehat{\theta}$ and $\Pi(l)$, respectively.
Based on these grounds, we would recommend using $\widehat{\theta}_{v}$ and $\Pi_{v}(l)$ in practice.

\section{Applications}

In this section, we consider two applications on the U.S. stock market.
Application 1 studies the low dimensional RCOV matrix series calculated by composite realized kernels (CRK) in Lunde, Shephard and Sheppard (2016).
Application 2 studies the high dimensional RCOV series calculated by TARVM estimator in Tao et al. (2011). %
\subsection{Application 1}
In this application, we revisit the RCOV matrix data of Hewlett-Packard Development Company, L.P. (HPQ), International Business Machines Corporation (IBM) and Microsoft Corporation (MSFT) in Lunde, Shephard and Sheppard (2016). This data set, denoted by $\{Y_t\}_{t=1}^{1474}$,
ranges from January 2006 to December 2011 with 1474  observations in total.
Here, two flash crashes are flagged in 6 May, 2010 and 9 August, 2011 and replaced by an average of the nearest five preceding and following matrices.

\begin{figure}[!htp]
	\centering
	\centerline{\includegraphics[scale=0.4]{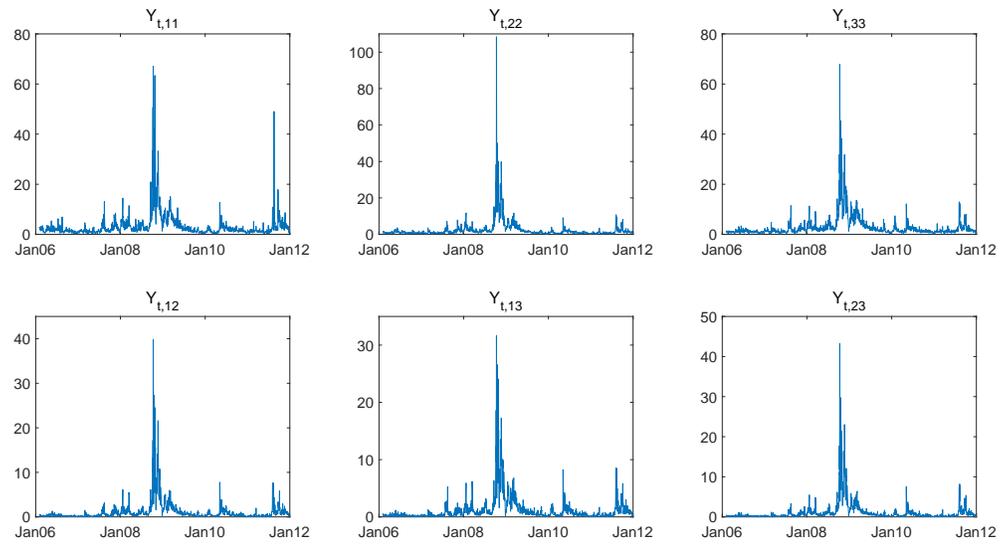}}
	\caption{\label{figapp1_1} Components of $Y_t$.}
\end{figure}

\begin{figure}[!htp]
	\centerline{\includegraphics[scale=0.4]{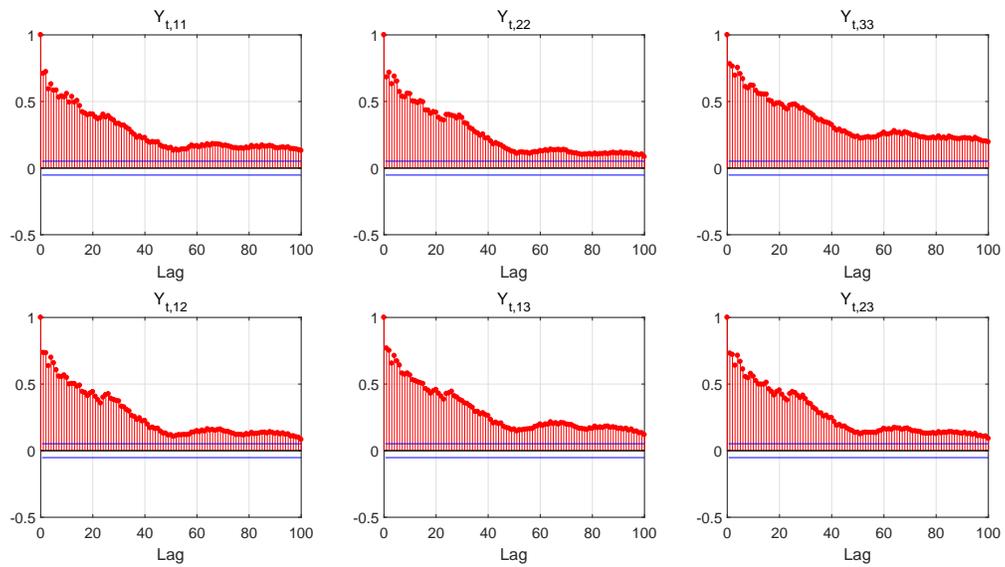}}
	\caption{\label{figapp1_2} Sample ACFs of each component $Y_{t,ij}$}
\end{figure}

Figure \ref{figapp1_1} plots the diagonal and off-diagonal components of $\{Y_{t}\}_{t=1}^{1474}$, exhibiting that
$Y_t$ has a clear clustering feature. Meanwhile, Figure \ref{figapp1_2} plots their sample autocorrelation functions (ACFs), which show the significant temporal dependence of $Y_t$. Based on these facts, we first fit $\{Y_{t}\}_{t=1}^{1474}$ by a diagonal VT-CBF model with $(P,Q,K)=(3,1,1)$, where
the order $K$ is taken as one for ease of model identification, and
the orders $P$ and $Q$ are selected by the Bayesian information criterion (BIC).
Specifically, this diagonal VT-CBF model is estimated using the two-step estimation procedure, and the corresponding estimates are
give in Table \ref{lowdimension}.
Second, since the sample ACFs of each component in Figure \ref{figapp1_2} decay slowly,
we also fit
$\{Y_{t}\}_{t=1}^{1474}$ by a diagonal VT-CBF-HAR model, and the related estimation results are also listed
in Table \ref{lowdimension}.
From this table, we find that the estimates of the degrees of freedom (especially for $\nu_2$) in both fitted models are close to each other, and
both estimates of $\nu_2$ are small indicating the heavy-tailedness of the examined data.
For the estimates of the mean parameter matrix $S$, its standard errors based on the VT-CBF model are smaller than those based on the
VT-CBF-HAR model. For other estimates of parameter matrices, the estimated diagonal components in each parameter matrix seem to have close values, meaning that the examined three stocks possibly have  similar temporal structures. This similarity can also be
seen from the values of persistence of each stock in Table \ref{lowdimension}, where the persistence of stock $s$
is defined by $\sum_{i=1}^{P}A_{1i,ss}^2+\sum_{j=1}^{Q}B_{1j,ss}^2$ for the VT-CBF model and
$A_{(d),ss}^2+A_{(w),ss}^2+A_{(m),ss}^2$ for the VT-CBF-HAR model. After estimation, we then apply our test statistics $\Pi_v(l)$  to both fitted models, and the results summarized in Table \ref{A1_innertest} imply that both fitted models are adequate at the 5\% level.

\begin{table}[htp]
	\centering
	\caption{The results of the estimated diagonal VT-CBF and VT-CBF-HAR models}
   \addtolength{\tabcolsep}{-1pt}
    \label{lowdimension}
	\begin{tabular}{ccllccccc}
		\hline
		\multicolumn{9}{c}{Diagonal VT-CBF model}                                                                                                                                                                    \\
		$\widehat{\nu}_{v}$      & \multicolumn{3}{c}{$\widehat{S}_{v}$}                              & $\widehat{A}_{11,v}$    & $\widehat{B}_{11,v}$    & $\widehat{B}_{12,v}$    & $\widehat{B}_{13,v}$    & persistence \\ \hline
		$74.0110$ & $3.1523$ & $1.1099$ & $1.1635$ & $0.7207$ & $0.5358$ & $0.0117$ & $0.4129$ & 0.9771      \\
        $(10.7545)$ & $(1.8844)$ & $(0.9031)$ & $(0.7705)$ & $(0.0223)$ & $(0.0365)$ & $(0.0176)$ & $(0.0354)$ &       \\
		$40.5849$  & $1.1099$ & $2.3683$ & $1.0965$ & $0.7200$ & $0.5620$ & $0.0119$ & $0.3800$ & 0.9788      \\
        $(3.9787)$  & $(0.9031)$ & $(2.1165)$ & $(0.9209)$ & $(0.0246)$ & $(0.0289)$ & $(0.0177)$ & $(0.0382)$ &      \\
		& $1.1635$ & $1.0965$ & $2.7883$ & $0.7118$ & $0.5579$ & $0.0127$ & $0.3977$ & 0.9762      \\
        & $(0.7705)$ & $(0.9209)$ & $(1.3276)$ & $(0.0211)$ & $(0.0292)$ & $(0.0190)$ & $(0.0354)$ &       \\ \hline \\
		\multicolumn{8}{c}{Diagonal VT-CBF-HAR model}                                                                                                                                                  &             \\
		$\widehat{\nu}_{v}$     & \multicolumn{3}{c}{$\widehat{S}_{v}$}                             & $\widehat{A}_{(d),v}$ & $\widehat{A}_{(w),v}$ & $\widehat{A}_{(m),v}$ & persistence         &             \\ \cline{1-8}
		$69.0222$  & $3.1523$ & $1.1099$ & $1.1635$ & $0.6954$ & $0.5735$ & $0.3891$ & 0.9639              &             \\
        $(6.2261)$  & $(2.2543)$ & $(1.0464)$ & $(0.8881)$ & $(0.0256)$ & $(0.0443)$ & $(0.0344)$ &               &             \\
		$40.4021$  & $1.1099$ & $2.3683$ & $1.0965$ & $0.6884$ & $0.6027$ & $0.3557$ & 0.9637              &             \\
        $(2.9408)$  & $(1.0464)$ & $(2.3391)$ & $(1.0210)$ & $(0.0275)$ & $(0.0318)$ & $(0.0426)$ &               &             \\
		& $1.1635$ & $1.0965$ & $2.7883$ & $0.6703$ & $0.6041$ & $0.3812$ & 0.9596              &             \\
        & $(0.8881)$ & $(1.0210)$ & $(1.4971)$ & $(0.0279)$ & $(0.0318)$ & $(0.0364)$ &              &             \\
\hline
	\end{tabular}
  \begin{tablenotes}
      \item[\dag] {\scriptsize The asymptotic standard errors are given in the parenthesis.}
  \end{tablenotes}
\end{table}

\begin{table}[!htp]
	\centering
	\caption{The results of $\Pi_{v}(l)$ for the diagonal VT-CBF and VT-CBF-HAR models}
	\label{A1_innertest}
	\begin{tabular}{llllll llllll}
		\hline
		&\multicolumn{5}{c}{Diagonal VT-CBF model} && \multicolumn{5}{c}{Diagonal VT-CBF-HAR model}\\
		\cline{2-6} \cline{8-12}
		$l$    & 2      & 3      & 4      & 5      & 6 &    & 2      & 3      & 4      & 5      & 6    \\
		\cline{2-6} \cline{8-12}
		$\Pi_v(l)$ & 1.494 & 4.170 & 8.004 & 9.428 & 11.513 &  & 4.385      &  6.127     & 7.004      &  10.310    &11.583 \\
		\cline{2-6} \cline{8-12}
			p-value & 0.474 & 0.244 & 0.091 & 0.093 & 0.074 &  & 0.112      &0.106      &0.136       & 0.067      &0.072 \\
		\hline
	\end{tabular}
\end{table}

Next, we consider the forecasting performance of our proposed diagonal VT-CBF and  VT-CBF-HAR models.
Specifically, we compute the 1-step, 5-step and 10-step predictions of the RCOV matrices, based on a rolling window procedure with window size equal to $T_0=800$.
That is, for $T_0\leq t\leq T-t_0$, we fit models based on $T_0$ observations $\{Y_s\}_{s=t-T_0+1}^{t}$, and forecast $\widehat{Y}_{t+t_0}$ with $t_0=1,5,10$ and calculate the forecasting error  by $\widehat{Y}_{t+t_0}-Y_{t+t_0}$.
To examine the importance of $\nu_2$ in the CBF models, we also apply the
diagonal VT-CAW and VT-CAW-HAR models to do prediction for the purpose of comparison.
The diagonal VT-CAW and VT-CAW-HAR models are defined in the same way as
the diagonal VT-CAW and VT-CAW-HAR models, except that the matrix-F distribution for $\Delta_t$ in the latter two models is replaced by
the Wishart distribution. Besides the CAW-type models, we further include a diagonal VAR-HAR model for comparison,
where this VAR model uses an HAR structure with
the diagonal autoregressive parameter matrices to fit $y_{t}=vech(Y_{t})$.

Table \ref{A1_forecast} gives the average of
forecasting errors in Frobenius and spectral norms for all models. From this table, we can find that regardless of the prediction horizon, the
diagonal VT-CBF-HAR model always has the smallest forecasting error in both norms.
Moreover, we apply the DM test (Diebold and Mariano, 1995) to examine whether the
diagonal VT-CBF-HAR model has a significant forecasting accuracy over other four competing models.
The corresponding testing results are given in  Table \ref{A1_forecast}, and they show that
the VT-CBF-HAR model is significantly better than its four competing models in terms of 5-step and 10-step forecasts.
For 1-step forecasts, the VT-CBF-HAR and VT-CBF model models have comparable forecasting accuracy, and
the VT-CBF-HAR model is significantly better than the remaining three models at level 10\%. Note that
the VAR-HAR model always performs worst in all examined cases, and this is probably because the VAR-HAR model brutally
disentangles the matrix-structure of the RCOV matrices, which may have some intrinsic and useful value for forecasts.

\begin{table}[]
		\caption{Forecasting errors based on different models and the related DM testing results}
	\label{A1_forecast}
\setlength{\tabcolsep}{3.6mm}{
\begin{tabular}{lllllll}
	\hline
	& \multicolumn{2}{l}{1-step} & \multicolumn{2}{l}{5-step} & \multicolumn{2}{l}{10-step} \\ \cline{2-7}
	Diagonal Model   & Frobenius     & Spectral      & Frobenius     & Spectral      & Frobenius      & Spectral      \\ \hline
	VT-CBF-HAR       & 1.5284        & 1.4607     & 1.9725        & 1.8850    & 2.2108         & 2.1091     \\
	VT-CBF           & 1.5349        & $1.4664$     & $1.9955$        & $1.9069^{\dag}$     & $2.2802^{*}$         & $2.1755^{*}$     \\
	VT-CAW-HAR & $1.5383^*$        & $1.4703^{*}$     & $2.0029^{*}$        & $1.9147^{*}$     & $2.2864^{\diamond}$        & $2.1813^{\diamond}$    \\
	VT-CAW           & $1.5390$        & $1.4699$     & $2.0253^{\diamond}$        & $1.9351^{\diamond}$     & $2.3364^{\diamond}$         & $2.2286^{\diamond}$     \\

	VAR-HAR          & $1.6472^{\diamond}$        & $1.5661^{\diamond}$     & $2.1700^{\diamond}$        & $2.0626^{\diamond}$     & $2.6088^{\diamond}$        & $2.4711^{\diamond}$    \\ \hline
\end{tabular}
  \begin{tablenotes}
      \item[1] {\scriptsize The DM test is used to compare the prediction accuracy between the diagonal VT-CBF-HAR and
the other four competing models. The result of the each competing model is marked with
``$\dag$'', ``$*$'' or ``$\diamond$'', if the DM test implies the Diagonal VT-CBF-HAR model gives significantly more accurate
predictions than this competing model at level 10\%, 5\% or 1\%, respectively.}
  \end{tablenotes}
  }
\end{table}

%
%

\subsection{Application 2}

In this section, we consider intra-day data of 112 stocks from four major sectors constituting S\&P 500 index: 31 stocks from financial sector, 31 stocks from industrial sector, 25 stocks from health care sector, and 25 stocks from consumer discretionary sector, see Table \ref{stockname}. All intra-day price data are downloaded from Wharton Research Data Services (WRDS) database, and they are taken from 1 July, 2009 to 30 December, 2016, including 1890 non-missing dates of trading data in total.
\begin{table}[h]
	\centering
	\label{stockname}
	\caption{Symbol of Stocks in Application 2 }
	\scalebox{1}{
		\begin{tabular}{ccccc}
			\hline
			Number & Financial Sector & Industrial Sector & Health Care Sector & Consumer Discretionary Sector \\ \hline
			1      & AFL              & BA                & A                  & AZO                           \\
			2      & AIG              & CAT               & ABC                & BBY                           \\
			3      & ALL              & FLR               & ABT                & BWA                           \\
			4      & AXP              & FLS               & AET                & CCL                           \\
			5      & BAC              & GD                & BAX                & GPC                           \\
			6      & BBT              & GE                & BDX                & GPS                           \\
			7      & BEN              & GWW               & BMY                & HD                            \\
			8      & BK               & HON               & BSX                & HRB                           \\
			9      & BLK              & IR                & CAH                & JWN                           \\
			10     & C                & ITW               & CI                 & KMX                           \\
			11     & CMA              & LLL               & CVS                & KSS                           \\
			12     & COF              & LMT               & HUM                & LEG                           \\
			13     & GS               & LUV               & JNJ                & LEN                           \\
			14     & HIG              & MAS               & LH                 & LOW                           \\
			15     & JPM              & MMM               & LLY                & MCD                           \\
			16     & KEY              & NOC               & MCK                & NKE                           \\
			17     & LNC              & NSC               & MDT                & NWL                           \\
			18     & MCO              & PH                & MRK                & PHM                           \\
			19     & MET              & PNR               & PFE                & RL                            \\
			20     & MMC              & PWR               & PKI                & TGT                           \\
			21     & MTB              & RHI               & SYK                & TIF                           \\
			22     & PFG              & ROK               & TMO                & TJX                           \\
			23     & PGR              & ROP               & UNH                & VFC                           \\
			24     & PNC              & RSG               & VAR                & WHR                           \\
			25     & PRU              & RTN               & WAT                & YUM                           \\
			26     & RF               & SNA               &                    &                               \\
			27     & STI              & SWK               &                    &                               \\
			28     & STT              & TXT               &                    &                               \\
			29     & TMK              & UNP               &                    &                               \\
			30     & USB              & UPS               &                    &                               \\
			31     & WFC              & UTX               &                    &                               \\ \hline
	\end{tabular}}
	\begin{tablenotes}
		\footnotesize
		\item[1]Note: Full names of selected stocks can be found in https://www.slickcharts.com/sp500
	\end{tablenotes}
\end{table}
 Based on 100 times log of the price data, the daily RCOV matrices $\{Y_{t}\}_{t=1}^{1890}$ are calculated by the TARVM method in  Tao et al. (2011) for each sector.

\begin{figure}[!h]
    \begin{center}
        \includegraphics[width=35pc,height=20pc]{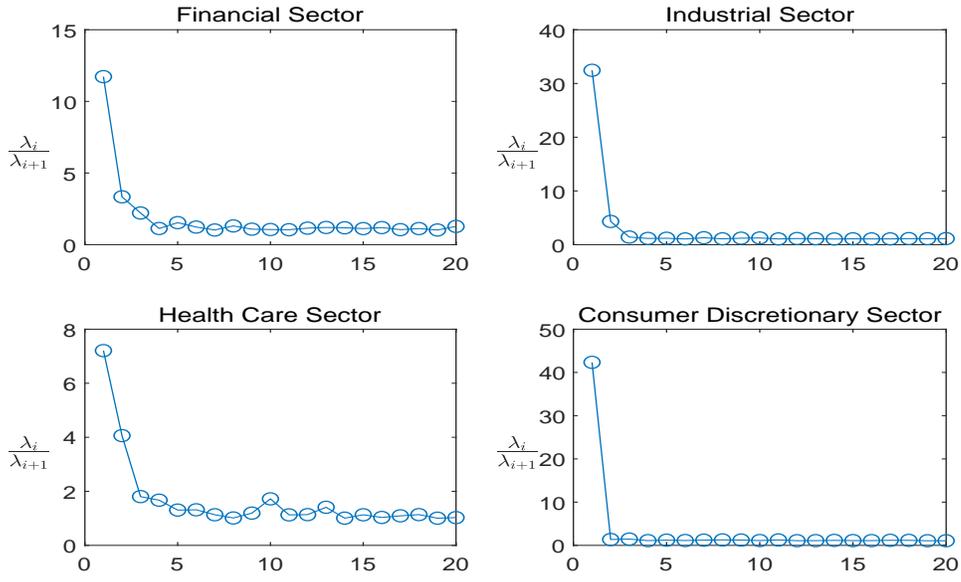}
    \end{center}
    \caption{\label{eigen} Ratios of adjacent eigenvalues of $\bar{S}$ for each sector}
\end{figure}

\begin{figure}[!h]
	\begin{center}
        \includegraphics[width=35pc,height=20pc]{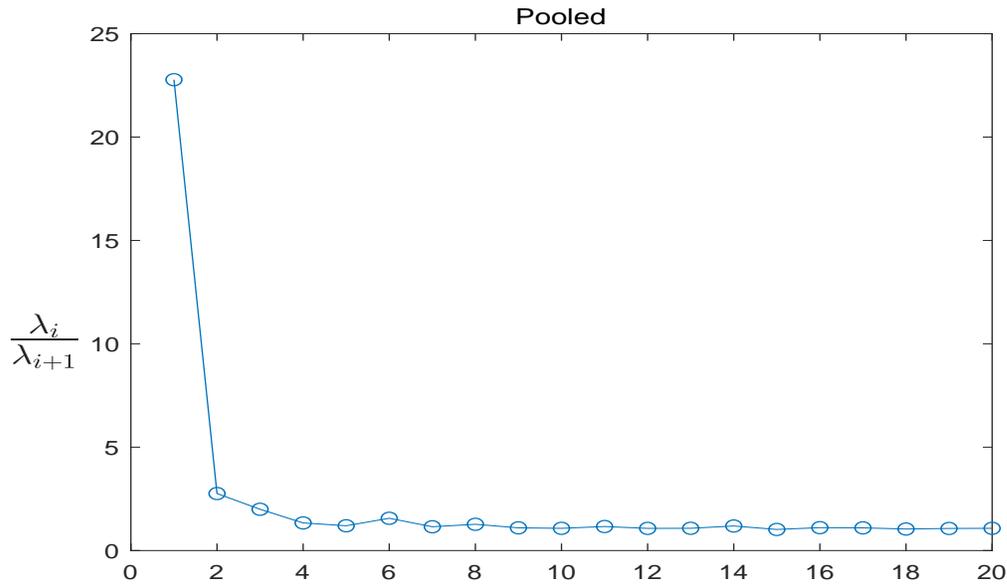}
	\end{center}
	\caption{\label{eigen_pool} Ratios of adjacent eigenvalues of $\bar{S}$ for all four pooled sectors}
\end{figure}

For each sector, since the dimension of the RCOV matrix is large, we fit the RCOV matrix data by the
diagonal F-VT-CBF and F-VT-CBF-HAR models. To do this,
we first look for the value of $r$ in model (\ref{factor_model}) by plotting the ratios $\{\frac{\lambda_{i}}{\lambda_{i+1}}\}$  for each sector in Fig\,\ref{eigen}, where $\{\lambda_{i}\}$ are the eigenvalues of $\bar{S}$ in descending order. From Fig\,\ref{eigen}, we can
choose $r=3$ for financial sector, $r=2$ for industrial sector, $r=2$ for health care sector, and $r=1$ for consumer discretionary sector.
To get more information, we also plot the ratios $\{\frac{\lambda_{i}}{\lambda_{i+1}}\}$ for all four pooled sectors
in Fig\,\ref{eigen_pool}, from which $r=3$ is suggested. This implies that
all 112 stocks considered may be driven by 3 latent factors, but among which only two may affect the industrial and health care sectors, and only one may affect the consumer discretionary sector. Hence, it is more reasonable to study the RCOV matrix data across sectors rather than together.

Next, we estimate the diagonal F-VT-CBF and F-VT-CBF-HAR models and choose the orders by a similar procedure as in Application 1, and the related
results are reported in Table \ref{est_results}. From this table, we can find that except for the mean parameter matrix,
the diagonal components of other parameter matrices seem to have different values, meaning that each component of  $Y_{ft}$
has a different dynamical structure. Moreover, the values of persistence for $Y_{ft,ss}$ show clear differences across four sectors, with
the largest persistence in financial sector and the smallest persistence in health care sector. This finding indicates that
the effect of past stock returns to its current volatility decays very slowly in the financial sector, while it
behaves oppositely in the health care sector.

\begin{table}[!h]
	\centering
	\caption{\label{est_results} The results of the estimated diagonal F-VT-CBF and F-VT-CBF-HAR models}
	\scriptsize\addtolength{\tabcolsep}{-4pt}
	\begin{tabular}{cccccccccccc}
		\hline
		                                                                            & \multicolumn{11}{c}{Diagonal F-VT-CBF model}                                                                                                                                                                                                                   \\
		Sector & $\widehat{\nu}_{fv}$    & \multicolumn{3}{c}{$\widehat{S}_{fv}$}                                 & $\widehat{A}_{11,fv}$        & $\widehat{B}_{11,fv}$    & $\widehat{B}_{12,fv}$   & $\widehat{B}_{13,fv}$   & $\widehat{B}_{14,fv}$   & persistence \\ \hline
		\multirow{6}{*}{Financial}                                                        & $35.3380$  & $25.7553$ & $0.6808$  & $0.1389$  & $0.7269$                        & $0.5118$  & $0.2741$ & $0.3219$ &                     & 0.9691      \\
		& ${(2.9679)}$  & ${(11.0314)}$ & ${(2.3577)}$  & ${(0.6519)}$  & ${(0.0348)}$                        & ${(0.0518)}$  & ${(0.1014)}$ & ${(0.0606)}$ &                     &     \\
		& $19.257$  & $0.6808$   & $2.5799$  & $0.0211$  & $0.6844$                        & $0.5382$  & $0.3172$ & $0.3628$ &                     & 0.9903      \\
				& ${(1.0419)}$  & ${(2.3577)}$   & ${(9.6931)}$  & ${(0.1730)}$  & ${(0.0608)}$                        & ${(0.1181)}$  & ${(0.1831)}$ & $(0.0699)$ &                     &    \\
		&                      & $0.1389$   & $0.0211$  & $1.6309$  & $0.7292$                        & $0.3010$  & $0.4490$ & $0.3817$ &                     & 0.9696      \\
			&                      & ${(0.6519)}$   & ${(0.1730)}$  & ${(1.8857)}$  & ${(0.0732)}$                        & ${(0.1468)}$  & ${(0.0897)}$ & ${(0.1201)}$ &                     &       \\	
		\hline
		\multirow{4}{*}{Industrial}                                                       & $24.9287$ & $17.3161$  & $2.1513$  &                      & $0.7277$                        & $0.6488$  &                     &                     &                     & 0.9505      \\
		 & ${(6.9460)}$ & ${(7.0877)}$  & ${(1.0290)}$  &                      & ${(0.0729)}$                        & ${(0.0709)}$  &                     &                     &                     &     \\
		& $22.7808$ & $2.1513$   & $1.0614$  &                      & $0.6716$                        & $0.6921$  &                     &                     &                     & 0.9300      \\
			& ${(7.6622)}$ & ${(1.0290)}$   & ${(0.3786)}$  &                      & ${(0.0317)}$                        & ${(0.0373)}$  &                     &                     &                     &       \\
		 \hline
		\multirow{4}{*}{Health Care}                                                      & $24.3415$ & $8.6744$   & $3.4402$  & & $0.7617$    & $0.5396$  & $0.1129$ &  &  & 0.8841       \\
		                            & ${(4.9720)}$ & ${(2.9442)}$   & ${(0.7505)}$  & & ${(0.1324)}$    & ${(0.0651)}$  & ${(0.6685)}$ &  &  &       \\
		& $ 15.9965$ & $3.4402$   & $2.185$  &   & $0.7351$    & $0.5706$  & $ 0.0001$ & &  & 0.8660          \\
		
			& $ {( 5.1757)}$ & ${(0.7505)}$   & ${(0.4998)}$  &   & ${(0.1407)}$    & ${(0.1585)}$  & $ {(0.8598)}$ & &  &          \\ \hline
		\multirow{4}{*}{\begin{tabular}[c]{@{}c@{}}Consumer\\ Discretionary\end{tabular}} & $22.4570$ & $15.3282$  &                      &                      & $0.7516$                        & $0.4517$  & $0.2604$ & $0.1971$ & $0.2666$ & 0.9467      \\
		& ${(4.0371)}$ & ${(4.9315)}$  &                      &                      & ${(0.0261)}$                        & ${(0.0724)}$  & ${(0.1171)}$ & ${(0.1711)}$ & ${(0.1032)}$ &       \\
		& $12.2757$ &                       &                      &                      &                                            &                      &                     &                     &                     &             \\
				& ${(1.4843)}$ &                       &                      &                      &                                            &                      &                     &                     &                     &             \\
		 \hline\\
		& \multicolumn{8}{c}{Diagonal F-VT-CBF-HAR model}                                                                                                                                                      &                     &                     &             \\
		Sector & $\widehat{\nu}_{fv}$   & \multicolumn{3}{c}{$\widehat{S}_{fv}$}                                & $\widehat{A}_{(d),fv}$ & $\widehat{A}_{(w),fv}$ & $\widehat{A}_{(m),fv}$ & persistence         &                     &                     &             \\ \cline{1-9}
		\multirow{6}{*}{Financial}                                                        & $38.0409$ & $25.7553$ & $0.6808$  & $0.1389$  & $0.7041$  & $0.5069$  & $0.4573$  & 0.9618              &                     &                     &             \\
		              & ${(3.1046)}$ & ${(15.5296)}$ & $(2.6814)$  & ${(0.8796)}$  & ${(0.0259)}$  & ${(0.0830)}$  & ${(0.1098)}$  &              &                     &                     &             \\
		& $18.9242$ & $0.6808$   & $2.5799$ & $0.0211$  & $0.6676$  & $0.4588$  & $0.5739$  & 0.9855              &                     &                     &             \\
		& ${(0.8746)}$ & ${(2.6814)}$   & ${(10.6104)}$ & ${(0.2816)}$  & ${(0.0441)}$  & ${(0.1162)}$  & ${(0.0628)}$  &              &                     &                     &             \\
		&                      & $0.1389$   & $0.0211$  & $1.6309$  & $0.7659$  & $0.2502$  & $0.5476$  & 0.9491              &                     &                     &             \\
				&                      & ${(0.8796)}$   & ${(0.2816)}$  & ${(1.2904)}$  & ${(0.0484)}$  & ${(0.1678)}$  & ${(0.0537)}$  &               &                     &                     &             \\
		 \cline{1-9}
		\multirow{4}{*}{Industrial}                                                       & $25.0002$ & $17.3161$ & $2.1513$  &                      & $0.7161$  & $0.5494$  & $0.3549$  & 0.9406              &                     &                     &             \\
		              & ${(5.9220)}$ & ${(10.0000)}$ & ${(1.2538)}$  &                      & ${(0.0699)}$  & ${(0.0758)}$  & ${(0.0458)}$  &              &                     &                     &             \\
		& $22.3305$ & $2.1513$   & $1.0614$  &                      & $0.6361$  & $0.6086$  & $0.3283$  & 0.8830              &                     &                     &             \\
				& ${(6.7511)}$ & ${(1.2538)}$   & ${(0.4310)}$  &                      & ${(0.0462)}$  & ${(0.0970)}$  & ${(0.1484)}$  &              &                     &                     &             \\
		 \cline{1-9}
		\multirow{4}{*}{Health Care}                                                      & $23.3766$ & $8.6744$   & $3.4402$  &  & $0.7259$  & $0.5357$  & $0.1944$  & 0.8625              &                     &                     &             \\
		 & ${(3.6648 )}$ & ${(3.2870)}$   & ${(0.8134)}$  &  & ${(0.1095)}$  & ${(0.1141)}$  & ${(0.0369)}$  &             &                     &                     &             \\
		& $ 16.1320$ & $3.4402$   & $2.1850$  &   & $0.6961$  & $0.5689$  & $0.0691$  & 0.8130               &                     &                     &             \\
				& $ {(4.6804)}$ & ${(0.8134)}$   & ${(0.5280)}$  &   & ${(0.0918)}$  & ${(0.1620)}$  & ${(0.2421)}$  &                &                     &                     &             \\
 \cline{1-9}
		\multirow{4}{*}{\begin{tabular}[c]{@{}c@{}}Consumer\\ Discretionary\end{tabular}} & $23.1216$ & $15.3282$  &                      &                      & $0.7285$  & $0.4865$  & $0.4092$  & 0.9348              &                     &                     &             \\
		& ${(3.2789)}$ & ${(6.0954)}$  &                      &                      & ${(0.0299)}$  & ${(0.0599)}$  & ${(0.0502)}$  &               &                     &                     &             \\
		& $11.9375$ &                       &                      &                      &                      &                      &                      &                     &                     &                     &             \\
				& ${(1.1630)}$ &                       &                      &                      &                      &                      &                      &                     &                     &                     &             \\ \hline
	\end{tabular}
  \begin{tablenotes}
      \item[\dag] {\scriptsize The asymptotic standard errors  given in the parenthesis are based on process $\widehat{Y}_{ft}$ rather than $Y_{ft}$.}
  \end{tablenotes}
\end{table}

In the end, we examine the forecasting performance of our F-CBF models. As in Application 1,
five different diagonal factor models (see Table \ref{pred_results}) are considered to forecast $Y_t$, based on a rolling window procedure with window size equal to $1000$. Their forecasting performance is evaluated by the average of forecasting errors in Frobenius  and spectral norms
as well as the results of the related DM test in Table \ref{pred_results}. From this table, we can see that  except for the health care sector,
the diagonal F-VT-CBF-HAR model always has the smallest forecasting error and
the diagonal F-VAR-HAR model has the largest forecasting error. For 1-step forecasts in the health care sector, the diagonal F-VT-CAW-HAR has slightly smaller forecasting error compared with the diagonal F-VT-CBF-HAR model.
In view of the results of DM test, the
 diagonal F-VT-CBF-HAR model has a significantly better performance than the other four competing models
in terms of 5-step and 10-step forecasts, but this advantage is slightly weak in terms of 1-step forecasts,
for which the diagonal F-VT-CBF and F-VT-CAW-HAR models have similar  performance in the industrial sector, and the
diagonal F-VT-CAW-HAR and F-VAR-HAR models have comparative performance in the health care sector.

\begin{table}[]
    \centering
	\caption{\label{pred_results} Forecasting errors based on different factor models and the related DM testing results}
    \scriptsize\addtolength{\tabcolsep}{-0.1pt}
	\begin{tabular}{clcccccc}
		\hline
		&                & \multicolumn{2}{c}{1-step}                                                                                      & \multicolumn{2}{c}{5-step}                                                                                      & \multicolumn{2}{c}{10-step}                                                                                                                                \\ \cline{3-8}
		Sector                                                                            & \multicolumn{1}{c}{Diagonal Model} & Frobenius                                              & Spectral                                                 & Frobenius                                               & Spectral                                                  & Frobenius                                               & Spectral                                                                                             \\ \hline
		\multirow{5}{*}{Financial}                                                        & F-VT-CBF-HAR   & 8.7701 & 7.9339 & 10.4581 & 9.7229 & 11.0221 & 10.3200                                          \\
		& F-VT-CBF       & $8.8116$                                        & $7.9824^{\dag}$                                        & $10.6677^{*}$                                          & $9.9315^{\diamond}$                                          & $11.3503^{*}$                                          & $10.6713^{\diamond}$                                                                                    \\
		& F-VT-CAW-HAR   & $8.7865$                                        & $7.9644^*$                                             & $10.5183$                                             & $9.8144^{\dag}$                                          & $11.1072$                                             & $10.4575$                                                                                       \\
		&  F-VT-CAW       & $8.8354^*$                                             & $8.0248^{\diamond}$                                          & $10.7097^{*}$                                          & $10.0151^{*}$                                         & $11.5030^{\diamond}$                                          & $10.8786^{\diamond}$                                                                                    \\
		& F-VAR-HAR        & $8.8878^*$                                             & $8.0662^{*}$                                           & $11.1055^{\diamond}$                                          & $10.4644^{\diamond}$                                         & $11.7725^{\diamond}$                                          & $11.1745^{\diamond}$                                                                                    \\ \hline
		\multirow{5}{*}{Industrial}                                                       &  F-VT-CBF-HAR   & 7.9567 & 7.0936 & 9.3154  & 8.5480 & 9.8270  & \begin{tabular}[c]{@{}c@{}}9.0842\end{tabular} \\
		&  F-VT-CBF       & $7.9735$                                               & $7.1169$                                               & $9.4094$                                              & $8.6334$                                             & $9.9837$                                           & $9.2397$                                                                                     \\
		& F-VT-CAW-HAR   & $7.9680$                                               & $7.1112^{\dag}$                                        & $9.4106^{*}$                                           & $8.6494^{*}$                                          & $10.0565^{\diamond}$                                          & $9.3255^{\diamond}$                                                                                     \\
		&  F-VT-CAW       & $7.9995^{*}$                                           & $7.1450^{*}$                                           & $9.4645^{*}$                                           & $8.7001^{*}$                                          & $10.1157^{*}$                                          & $9.3826^{*}$                                                                                     \\
		&  F-VAR-HAR        & $8.0567^*$                                             & $7.2170^*$                                             & $9.6801^{\diamond}$                                           & $8.9531^{\diamond}$                                          & $10.2809^{\diamond}$                                          & $9.5794^{\diamond}$                                                                                     \\ \hline
		\multirow{5}{*}{Health Care}                                                      &  F-VT-CBF-HAR   & $6.6253$                                               & $5.8586$                                               & 7.4977                 & 6.8076                & 7.8436                  & 7.1863                                                            \\
		&  F-VT-CBF       & $6.6628^{\dag}$                                           & $5.9019^{\dag}$                                             & $7.6400^{*}$                                           & $6.9605^{*}$                                          & $8.0708^{\diamond}$                                           & $7.4398^{\diamond}$                                                                                     \\
		&  F-VT-CAW-HAR   & 6.6126                 & 5.8559                 & $ 7.5658^{*}$                                           & $ 6.8892^{*}$                                           & $7.9743{\diamond}$                                           & $7.3317^{\diamond}$                                                                                     \\
		&  F-VT-CAW       & $6.7451^{\diamond}$                                          & $6.0117^{\diamond}$                                          & $8.0423^{\diamond}$                                           & $7.3944^{\diamond}$                                          & $ 8.3738^{\diamond}$                                           & $7.7569^{\diamond}$                                                                                     \\
		&  F-VAR-HAR        & $6.6688$                                               & $5.8954$                                               & $7.6163^{*}$                                           & $6.9389^{*}$                                          & $7.9457^{\dag}$                                           & $7.2872$                                                                                     \\ \hline
		\multirow{5}{*}{\begin{tabular}[c]{@{}c@{}}Consumer\\ Discretionary\end{tabular}} &  F-VT-CBF-HAR   & 8.3355                 & 7.0130                & 9.3278                  & 8.1225                 & 9.6830                  & 8.5081                                                            \\
		&  F-VT-CBF       & $8.3552^{\dag}$                                        & $7.0415^*$                                             & $9.4191^{\dag}$                                           & $8.2195^{\dag}$                                          & $9.8426^{*}$                                           & $8.6883^{*}$                                                                                     \\
		&  F-VT-CAW-HAR   & $8.3517^{*}$                                           & $7.0307^*$                                             & $9.3886^{\diamond}$                                           & $8.1935^{\diamond}$                                          & $9.7918^{\diamond}$                                           & $8.6294^{\diamond}$                                                                                     \\
		&  F-VT-CAW       & $8.3727^*$                                             & $7.0560^{\diamond}$                                          & $9.4489^{*}$                                           & $8.2546^{*}$                                          & $9.9211^{\diamond}$                                           & $8.7754^{\diamond}$                                                                                     \\
		&  F-VAR-HAR        & $8.3914^*$                                             & $7.0762^*$                                             & $9.5017^{*}$                                           & $8.3282^{*}$                                          & $9.9085^{\diamond}$                                           & $8.7575^{\diamond}$                                                                                     \\ \hline
	\end{tabular}
  \begin{tablenotes}
      \item[1] {\scriptsize The DM test is used to compare the prediction accuracy between the diagonal F-VT-CBF-HAR and
the other four competing models. The result of the each competing model is marked with
``$\dag$'', ``$*$'' or ``$\diamond$'', if the DM test implies the Diagonal F-VT-CBF-HAR model gives significantly more accurate
predictions than this competing model at level 10\%, 5\% or 1\%, respectively.}
  \end{tablenotes}
\end{table}

\section{Concluding Remarks}

This paper proposes a new CBF model to study the dynamics of the RCOV matrix.
For this CBF model, we explore its stationarity and moment properties, establish the asymptotics of
its maximum likelihood estimator, and investigate the inner-product-based tests for its model checking.
Hence, a systematic inferential tool of this CBF model is available for empirical researchers.
In order to deal with large dimensional RCOV matrices, we also construct two reduced CBF models:
the VT-CBF model and the F-CBF model.
For both reduced models, the asymptotic theory of the estimated parameters is derived.
Compared with the CAW model with Wishart innovations, the CBF model with matrix-F innovations is more able in capturing the heavy-tailed RCOV.
This advantage is demonstrated by two real examples on U.S. stock markets.
As motivated by Chiriac and Voev (2011), one obvious future work is to introduce the fractional integration structure into our CBF models.
Another interesting potential future work could extend
the idea of using the matrix-F innovation in
a number of ways resulting in a large family of models, which shall be important to
study the positive definite dynamics.

 \section*{Supplementary Material}
 The online Supplementary Material contains the proofs of all theorems, and some useful derivatives.


\end{document}